\documentclass[aap,nameyear,MSNbibl,dvips]{arximspdf}
\usepackage{mathrsfs}
\usepackage{graphics}
\doi{10.1214/09-AAP645}
\volume{20}
\issue{3}
\pubyear{2010}
\firstpage{978}
\lastpage{1004}

\makeatletter

\newproclaim{remark}{Remark}[section]

\newtheorem{proposition}[remark]{Proposition}
\newtheorem{lemma}[remark]{Lemma}
\newtheorem{theorem}[remark]{Theorem}

\makeatother

\begin{document}
\begin{frontmatter}

\title{Asymptotic behavior of the rate of adaptation}
\runtitle{Rate of adaptation}

\begin{aug}
\author[A]{\fnms{Feng} \snm{Yu}\corref{}\thanksref{T1}\ead[label=e1]{feng.yu@bristol.ac.uk}},
\author[B]{\fnms{Alison} \snm{Etheridge}\ead[label=e2]{etheridg@stats.ox.ac.uk}} and
\author[C]{\fnms{Charles} \snm{Cuthbertson}\thanksref{T2}\ead[label=e3]{Charles.W.Cuthbertson@MorganStanley.com}\ead[label=u1,url]{www.foo.com}}
\runauthor{F. Yu, A. Etheridge and C. Cuthbertson}
\affiliation{University of Bristol, University of Oxford and Morgan Stanley}
\address[A]{F. Yu\\
School of Mathematics\\
University of Bristol\\
Bristol, BS8 1TW\\
United Kingdom\\
\printead{e1}} 
\address[B]{A. Etheridge\\
Department of Statistics\\
University of Oxford\\
1 South Parks Road\\
Oxford, OX1 3TG\\
United Kingdom\\
%
\printead{e2}}
\address[C]{C. Cuthbertson\\
Institutional Equity Division\\
Morgan Stanley\\
1565 Broadway\\
New York, New York 10036-8293\\
USA\\
\printead{e3}\\
\printead{u1}}
\end{aug}

\thankstext{T1}{Supported by EPSRC/GR/T19537 while at the University of
Oxford.}

\thankstext{T2}{Supported by EPSRC DTA.}

\received{\smonth{8} \syear{2007}}
\revised{\smonth{9} \syear{2009}}

%
\begin{abstract}
We consider the accumulation of beneficial and deleterious mutations in
large asexual populations. The rate of adaptation is affected by the
total mutation rate, proportion of beneficial mutations and population
size $N$. We show that regardless of mutation rates, as long as the
proportion of beneficial mutations is strictly positive, the adaptation
rate is at least $\mathcal{O}(\log^{1-\delta} N)$ where $\delta$ can be
any small positive number, if the population size is sufficiently
large. This shows that if the genome is modeled as continuous, there is
no limit to natural selection, that is, the rate of adaptation grows in
$N$ without bound.
\end{abstract}

%

\end{frontmatter}

\section{Background and introduction}
\label{sec1} We consider the accumulation of mutations in large asexual
populations. The mutations that biological organisms accumulate over
time can be classified into three categories: beneficial, neutral and
deleterious. Beneficial mutations increase the fitness of the
individual carrying the mutation, while deleterious mutations decrease
fitness; neutral mutations have no effect on fitness. Adaptation is
driven by accumulation of beneficial mutations, but it is limited by
\textit{clonal interference}
(clones that carry different beneficial
mutations compete with each other and interfere with the other's
growth in the population).
Fisher and Muller argued
for the importance of this effect as early as the 1930s [Fisher
(\citeyear{Fisher1930}),
Muller (\citeyear{Muller64})]. Here, we are concerned with the rate of
adaptation, that is
the rate of increase of mean fitness in the population.

The simplest scenario one can consider is one in which a single
beneficial mutation arises in an otherwise neutral population and no
further mutations occur until the fate of that mutant is known. This
situation is well understood. The most basic question one can ask is
what is $p_{\mathrm{fix}}$, the fixation probability of the mutation. This was
settled by Haldane (\citeyear{Haldane27}), who showed that under a
discrete generation haploid model, if the selection coefficient
associated with the mutation is $s$, then under these circumstances
$p_{\mathrm{fix}}\approx2s$. In this case, $p_{\mathrm{fix}}$ is almost independent of
the population size, $N$.

When the mutation does fix, the process whereby it increases in
frequency from $1/N$ to $1$ is known as a \textit{selective sweep}. The
duration of a selective sweep is $\mathcal{O}(\log(s N)/s)$
generations. If one assumes that the mutation rate per individual per
generation is $\mu$, then the overall mutation rate will be
proportional to population size and we see that for large populations
the assumption that no new mutation will arise during the timecourse of
the sweep breaks down. Instead one expects multiple overlapping sweeps.
In an asexual population, mutations can only be combined if they occur
sequentially within the same lineage. This means that, on the one hand,
alleles occurring on the same lineage can boost one another's chance of
fixation, but on the other hand alleles occurring on distinct lineages
competitively exclude one another. The net effect is to slow down the
progress of natural selection. This is an extreme form of the
Hill--Robertson effect. Hill and Robertson (\citeyear{HillRobertson66})
were the first to quantify the way in which linkage between two sites
under selection in a finite population (whether sexually or asexually
reproducing) limits the efficacy of natural selection. In a sexually
reproducing population, recombination breaks down associations between
loci and so ameliorates the Hill--Robertson effect, suggesting an
indirect selective force in favor of recombination. Further
quantitative analysis of the interference between selected loci is
provided by Barton (\citeyear{Barton95}) who considers the probability
of fixation of two favorable alleles in a sexually reproducing
population. His method is only valid if the selection coefficient of
the first beneficial mutation to arise is larger than that of the
second. Cuthbertson, Etheridge and Yu (\citeyear{YuEtCul08}) consider
the same question in the general setting. The conclusion from both
works is that fixation probabilities are reduced, sometimes
drastically, because of interference between the two mutations.
Furthermore, if the second mutation is stronger than the first, then
Cuthbertson, Etheridge and Yu (\citeyear {YuEtCul08}) show that the
strength of interference can be strongly dependent on population size.
In this work, we do not consider the effects of recombination, since we
only work with asexual populations.

Since all beneficial mutations eventually become either extinct or
ubiquitous in the population, the rate of adaptation, defined to be the
rate of increase of the mean fitness of the population, is proportional
to $\mu s p_{\mathrm{fix}} N$, where $\mu N$ is the total number of beneficial
mutations that occur to all individuals in the population in a single
generation and we assume $p_{\mathrm{fix}}$ to be the same for all beneficial
mutations, which is the case for the system in stationarity. If
$p_{\mathrm{fix}}$ is independent of population size, then we expect an
adaptation rate of $\mathcal{O}(N)$. However, a explained above, the
occurrence of simultaneous selective sweeps reduces $p_{\mathrm{fix}}$ and so
$p_{\mathrm{fix}}$ may not be $\mathcal{O}(1)$. This leads to the following
question: if one does not limit the number of simultaneous selective
sweeps, what is $p_{\mathrm{fix}}$, or equivalently, what is the rate of
adaptation? As $N\to\infty$, is the rate of adaptation finite or does
it increase without bound? There has been some controversy surrounding
this question. Some work [e.g., Barton and Coe
(\citeyear{BartonCoe07})] suggests that there is an asymptotic limit to
the rate of adaptation. Other authors [e.g., Rouzine, Wakeley and
Coffin (\citeyear{RouzineWC03}), Wilke (\citeyear {Wilke04}) and Desai
and Fisher (\citeyear{DesaiFisher07})] argue that no such limit exists.
Here, we study this problem in a mathematically rigorous
framework.

Previous work on this question has adopted two general approaches:
(i) calculate the fixation probability $p_{\mathrm{fix}}$ directly, and (ii)
study the distribution of fitness of all individuals in the
population and asks how this distribution evolves with time. The
first approach was used in Gerrish and Lenski (\citeyear
{GerrishLenski98}), Wilke (\citeyear{Wilke04})
and Barton and Coe (\citeyear{BartonCoe07}).
Gerrish and Lenski (\citeyear{GerrishLenski98}) were the first to
present a quantitative
analysis of
the rate of adaptation in the presence of clonal interference. They
obtained approximate integral expressions for the fixation
probability of a beneficial mutation and thus the expected rate of
adaptation. Orr (\citeyear{Orr00}) generalized the results of
Gerrish and Lenski (\citeyear{GerrishLenski98}) to include the effects
of deleterious
mutations. Wilke (\citeyear{Wilke04}) combined the works of Gerrish and
Lenski (\citeyear{GerrishLenski98}) and Orr (\citeyear{Orr00}) to
obtain approximate expressions for
the adaptation rate that grow logarithmically or doubly
logarithmically for large $N$. In all three works, the authors used
a sequence of approximations before arriving at an expression for
the fixation probability or the adaptation rate. It seems to be
highly nontrivial to turn any of these approximation steps into a
rigorous mathematical argument and so we do not follow their
approaches here.

The second approach, to consider the distribution of fitness in the
population, was used in Rouzine, Wakeley and Coffin (\citeyear
{RouzineWC03}), Brunet et al. (\citeyear{Brunetal06}) and Rouzine,
Brunet and Wilke (\citeyear{RouzineBW07}). As in the work described in
the last paragraph, Rouzine, Wakeley and Coffin (\citeyear
{RouzineWC03}) take fitness effects to be additive, but whereas before
the selection coefficient of each new mutation was chosen from a
probability distribution, now all selection coefficients are taken to
be equal. In this setting, a~beneficial and a deleterious mutation
carried by the same individual cancel one another out and an
individual's fitness can be characterized by the \textit{net} number of
beneficial mutations which it carries (which may be negative). Writing
$P_k$ for the proportion of individuals with fitness equivalent to $k$
beneficial mutations, $\{P_k\}_{k\in{\mathbb Z}}$ forms a type of
\textit{traveling wave} whose shape remains basically unchanged over
time. The position of the wave moves to the left or the right on the
fitness axis, depending on whether the adaptation rate is positive or
negative. This is similar to traveling waves arising from
reaction--diffusion equations in the PDE literature [see, e.g., Chapter
15 of Taylor (\citeyear{Taylor96})]. In the current setting, however,
the shape of the wave actually fluctuates stochastically even after a
long time. So the wave can be regarded as a stochastic traveling wave,
and its \textit{speed} is proportional to the rate of adaptation.
Rouzine, Wakeley and Coffin (\citeyear{RouzineWC03}) studied a
multilocus model that does not include recombination but does include
beneficial, deleterious and compensating mutations. They found that
the rate of adaptation (i.e., the speed of the traveling wave)
asymptotically depends logarithmically on population size~$N$, which is
consistent with results of in vitro studies of a type of RNA virus in
Novella et al. (\citeyear{Novellaetal95}, \citeyear{Novellaetal99}).
Rouzine, Brunet and Wilke (\citeyear{RouzineBW07}) presents the same
approach but with more detailed derivations and improved treatments of
the stochastic edge.

Desai and Fisher (\citeyear{DesaiFisher07}) also adopts the traveling
wave approach. Their method of studying the adaptation rate, however,
differs from that of Rouzine, Wakeley and Coffin
(\citeyear{RouzineWC03}) and Rouzine, Brunet and Wilke
(\citeyear{RouzineBW07}) in that they consider the fitness variation of
the population to be in mutation--selection balance, and ask how much
variance in fitness can the population maintain while this variation is
being selected on. The conclusion they reach is that this variation
(hence the adaptation rate) increases logarithmically with both
population size and mutation rate.

Brunet et al. (\citeyear{Brunetal06}) study a model in which each of
the $N$ individuals in the population gives birth to $k$ offspring,
each of which has a fitness that differs from the fitness of its parent
by a random amount and finally the $N$ fittest individuals are used to
form the next generation. This model resembles artificial selection,
rather than natural selection, but it may be easier to study because
the density of individuals of a certain fitness in the next generation
has a kind of \textit{local dependence} on that density in the current
generation. This is quite different from the behavior considered in
Rouzine, Wakeley and Coffin (\citeyear{RouzineWC03}) and our work in
this article, where the density of individuals of a given fitness
depends on the whole fitness distribution of the parental population.

This work originally arose from discussions with Nick Barton and
Jonathan Coe which focused on limits to the rate of adaptation when all
mutations are beneficial. In reality, most mutations are either neutral
or deleterious. In particular, if all mutations in an asexual
population were deleterious, then the population would irreversibly
accumulate deleterious mutations, a process known as Muller's ratchet.
The first mathematically rigorous analysis of Muller's ratchet is due
to Haigh (\citeyear{Haigh78}). There, a Wright--Fisher model is
formulated that incorporates the effects of selection and mutation.
Again all mutations carry equal weight so that individuals can be
classified according to how many mutations they carry. Haigh
(\citeyear{Haigh78}) showed that if the population size is infinite (so
that the dynamics of the model become deterministic) then there is a
stationary distribution. In the finite population case, however, this
is not the case. At any given time, there is a \textit{fittest class},
corresponding to those individuals carrying the smallest number of
mutations, but this class will eventually be lost due to genetic drift
(the randomness in the reproduction mechanism). This loss is permanent
since there is no beneficial or back mutation to create a class fitter
than the current fittest class. The next fittest class then becomes the
fittest class, but that will be lost eventually as well and the entire
population grows inexorably less fit. Higgs and Woodcock
(\citeyear{HiggsWoodcock95}) derived a set of moment equations for
Haigh's model but these are not closed and so are hard to analyse.
Instead, their results rely mainly on simulations. Stephan, Chao and
Smale (\citeyear{StephanChaoSmale93}) and Gordo and Charlesworth
(\citeyear{GordoCharlesworth00}) use (slightly different)
one-dimensional diffusions to approximate the size of the fittest
class. Etheridge, Pfaffelhuber and Wakolbinger (\citeyear{EPW07}) go
much further along this line (and provide a more thorough review of the
literature than that included here). They conjecture and provide
justification for a phase transition and power law behavior in the rate
of the ratchet. But in spite of the very considerable body of work on
Muller's ratchet, even a rigorous expression for the rate of decline in
mean fitness of the population remains elusive.

Muller's ratchet caricatures the evolution of a population in which
there is no recombination and no beneficial mutation. Such a
population is doomed to become progressively less and less fit. So
how can a species overcome the ratchet? If it reproduces sexually,
then recombination of parental chromosomes can create offspring that
are fitter than either parent and so Muller's ratchet has been
proposed as an explanation for the evolution of sexual reproduction
[e.g., Muller (\citeyear{Muller64}), Felsenstein (\citeyear{Felsenstein74})].
But not
all populations reproduce sexually. Another mechanism which has the
potential to overcome Muller's ratchet is the presence of beneficial
mutations, and it is this mechanism that we shall consider here.
More specifically, we pose the following question: with both
beneficial and deleterious mutations, does a sufficiently large
population overcome Muller's ratchet?

The conclusion we reach, through both nonrigorous (Section \ref{sec3})
and rigorous (Theorem \ref{thm:main}) approaches, is the following:
as long as the proportion of beneficial mutations is strictly
positive, the rate of adaptation is roughly $\mathcal{O}(\log N)$ for
large~$N$, where $N$ is the population size and
time is measured in generations. This shows that even
with a tiny proportion of beneficial mutations, a large enough
population size will yield a positive adaptation rate, in which case
the entire population grows fitter at a high rate and Muller's
ratchet is overcome.
It also shows, in particular, that the rate of adaptation grows without bound
as $N\to\infty$ in the all-mutations-beneficial case. This is
consistent with the findings of Rouzine, Wakeley and Coffin (\citeyear
{RouzineWC03}) and
Wilke (\citeyear{Wilke04}).

\begin{figure}

\includegraphics{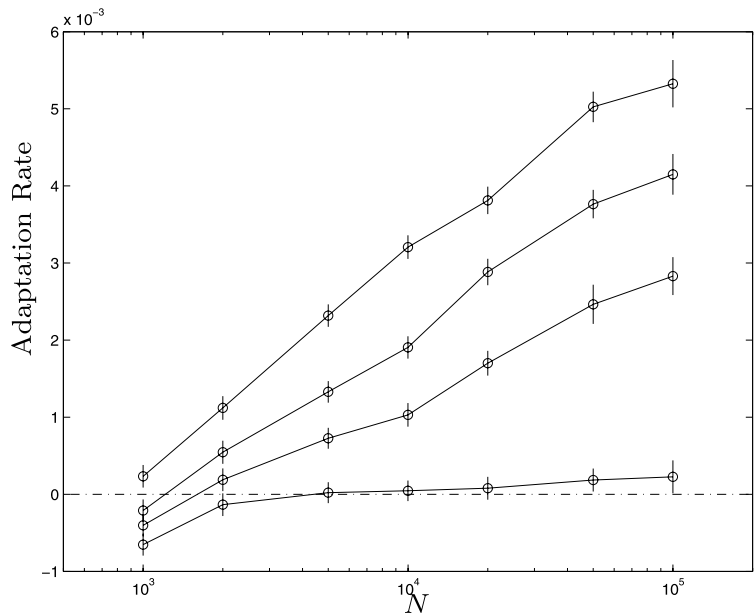}

\caption{Adaptation rate against population size, from top to bottom,
for $q=4\%, 2\%, 1\%$ and $0.2\%$, $\mu= 0.01$ and
$s= 0.01$.
Circles represent data points obtained from simulation,
$q$ is the probability that a mutation is advantageous
and vertical bars represent one standard deviation.}
\label{fig:means}
\end{figure}

Figure \ref{fig:means} plots the adaptation rate against log
population size
from simulation results of the model we consider in this article.
We observe that for each set of parameters $q$, $\mu$ and $s$,
the rate of adaptation is roughly proportional to $\log N$ and small
population sizes may result in negative adaptation rates.
Furthermore, larger $q$ results in a higher adaptation rate for
fixed $\mu$ and $s$. The upshot is that with $\mu$ and $s$
held constant, a smaller proportion of beneficial mutations needs
a larger population size for Muller's ratchet to be overcome.


In the model, we study here the selection coefficient $s$ is held fixed as
$N\rightarrow\infty$. This is known as a ``strong selection'' model.
Our interest is in the behavior of the model for very large $N$. It is
not clear in this setting how to pass to an infinite population limit and
so we must work with a model based on discrete individuals. An alternative
model, the so-called weak selection model, is used to address behavior of
very large populations when $Ns$ is not too large. By fixing $Ns$ (as
opposed to $s$), one can pass to an infinite population limit. The limiting
model comprises a countably infinite system of coupled stochastic
differential equations for the frequencies of individuals of different
fitnesses within the population. Preliminary calculations for this model
are presented in Yu and Etheridge (\citeyear{YuEtheridge08}).

%

This work is organized as follows. In Section \ref{sec2}, we formulate
our model. In the biological literature, one would expect to see a
Wright--Fisher model but since we are interested in large populations,
we expect the same results for the much more mathematically tractable
Moran particle model which we describe. We also perform some
preliminary calculations. In Section \ref{sec3}, we present a
nonrigorous argument that leads to an asymptotic adaptation rate
of roughly $\mathcal{O}(\log N)$. In Section~\ref{sec4}, we present and
prove our main rigorous result that establishes a lower bound of
$\log^{1-\delta} N$ for any $\delta>0$ on the adaptation rate. And
finally, in Section \ref{sec5}, we prove the supporting lemmas required
for the proof of our main theorem.

\section{The finite population Moran model}
\label{sec2}
We assume constant population size~$N$. For each $N \in\mathbb{N}$,
let $X_i(t) \in{\mathbb Z}$, $i=1,\ldots,N$, denote the \textit{fitness type}
of the $i$th individual, defined to be the
number of beneficial mutations minus the number of deleterious
mutations carried by the individual.
For $k\in{\mathbb Z}$, let $P_k(t)$ denote the
proportion of individuals that have fitness type $k$ at time $t$, that is,
\[
P_k=\frac1 N \sum_{i=1}^N \mathbf{1}_{\{ X_i = k \}}.
\]
We use $\mathscr{P}^{(N)}({\mathbb Z})$ to denote the space of probability
measures $p$ on ${\mathbb Z}$ formed by $N$ point masses each with
weight $1/N$,
and define
\[
S^{(N)} = \mathscr{P}^{(N)}({\mathbb Z})
\]
to be the state space for $P_k(t)$.
For $p\in S^{(N)}$, we define $p_k=p(\{k\})$ and
%
%
\begin{eqnarray} \label{def:cn}
p_{[k,l]} &=& \sum_{i=k}^l p_i, \nonumber\\
m(p) &=& \langle k,p\rangle= \sum_{k\in{\mathbb Z}} k p_k,
\\
c_n(p) &=& \sum_{k\in{\mathbb Z}} \bigl(k-m(p)\bigr)^n p_k.\nonumber
\end{eqnarray}
In particular, $m(p)$ is the mean fitness of the population, and
$c_2(p) = \langle k^2,p\rangle- \langle k,p\rangle^2$ is the
$2$nd
central moment of the population fitness, that is, its variance.
We sometimes abuse notation and use $P$ to denote the probability mass function
of different fitness types associated with the probability measure $P$.

The model of interest is one where each individual accumulates
beneficial mutations at a Poisson rate $q \mu$ and deleterious
mutations at rate $(1-q) \mu$. We assume a so-called infinitely-many-loci
model where each mutation is assumed to be new and occur at a different
locus on the genome.
All individuals experience selection
effects via a selection mechanism
(which introduces a drift reflecting the
differential reproductive success based on fitness)
and the effect of genetic drift
via a resampling mechanism. The mechanisms of this model are described below:
\begin{enumerate}
\item \textit{Mutation}: for each individual $i$, a mutation event
occurs at rate $\mu$.
With probability $1 - q$, $X_i$ changes to $X_i - 1$ and with
probability $q$, $X_i$ changes to $X_i +1$.
\item \textit{Selection}: for each pair of individuals $(i, j)$, at rate
$\frac{s}{N}(X_i - X_j)^+$, individual $i$ replaces individual $j$.
\item \textit{Resampling}: for each pair of individuals $(i, j)$, at
rate $\frac{1}{N}$, individual $i$ replaces individual $j$.
\end{enumerate}
This model has a time scale such that one unit of time corresponds
roughly to one generation. A more sophisticated model should
consider mutations that have a distribution of fitness effects, for example,
an independent exponentially distributed selective advantage
associated with each new beneficial mutation as proposed by
Gillespie (\citeyear{Gillespie91}).
Recent work by Hegreness et al. (\citeyear{Hegrenessetal06}), however,
suggests that in
models where beneficial mutations have a distribution of fitness
advantages, evolutionary dynamics, for example, the distribution of successful
mutations which ultimately determines the rate of adaptation,
can be reasonably described by an
equivalent model where all beneficial mutations confer the same
fitness advantage.
One can also describe the mechanisms in the above model in terms of
the~$P_k$'s,
\begin{enumerate}
\item \textit{Mutation}: for any $k\in{\mathbb Z}$,
at rate $(1-q)\mu N P_k$, $P_k$ decreases by $\frac1 N$
and $P_{k-1}$ increases by $\frac1 N$;
at rate $q\mu N P_k$, $P_k$ decreases by $\frac1 N$
and $P_{k+1}$ increases by $\frac1 N$.
\item \textit{Selection}: for any pair of $k,l\in{\mathbb Z}$ such
that $k>l$,
at rate $s(k - l) N P_k P_l$, $P_k$ increases by $\frac1 N$ and
$P_l$ decreases by $\frac1 N$.
\item \textit{Resampling}: for any pair of $k,l\in{\mathbb Z}$, at rate
$N P_k P_l$, $P_k$ increases by $\frac1 N$
and $P_l$ decreases by $\frac1 N$.
\end{enumerate}
We use $(P,X)$ to denote the process evolving under the above mechanism,
where $X$ describes the fitness types of the
$N$ exchangeable individuals and $P$ describes the empirical measure formed
by the fitness types of these individuals.
If there is no confusion, we drop $X$ and denote the process simply by $P$.
The main result of this work, Theorem \ref{thm:main}, states that under
the above model, the mean fitness increases at a rate of at least
$\mathcal{O}(\log^{1-\delta} N)$ for any $\delta>0$ after a sufficiently
long time.\looseness=-1
\begin{remark}\label{rem21}
Notice that the resampling acts on \textit{ordered} pairs, so that the overall
rate at which an individual is affected by a resampling event is $2N$ and
at such an event it has equal chance of reproducing or dying. It would be
more usual to have resampling at half this rate, but this choice of
timescale does not change the results and
will save us many factors of two later.

Often one combines the resampling and selection into a single term.
Each pair of individuals is involved in a reproduction event at some constant
rate and the effect of selection is then that it is more likely to
be the fitter individual that reproduces. Since $s$ is typically
rather small, our simpler formulation
is a very small perturbation of this model and again the statement of
our results would not be changed in that framework.
\end{remark}
\begin{remark}
We take the selection mechanism to be additive instead of multiplicative,
that is, the fitness type of an individual with $k$ beneficial mutations is
$1+s k$ instead of $(1 + s)^k$. Even though
$(1 + s)^k \approx1 + s k$ is only valid for small $s$ and $k$,
$(1 + s)^k \ge1 + s k$ holds for all $s\in[-1,\infty)$, thus
our main result of a lower bound on the rate of adaptation
also holds for multiplicative selection effects.
\end{remark}

One can construct the process $X(t)$ using Poisson random measures and
Poisson processes. More specifically, let $l$ denote the Lebesgue measure
on ${\mathbb R}$. For each $i\in{\mathbb Z}$,
let $\Lambda^{(N)}_{1i,b}$ and $\Lambda^{(N)}_{1i,d}$ be independent
Poisson processes with intensities $q\mu$ and $(1-q)\mu$, respectively.
For each $i,j\in{\mathbb Z}$, let $\Lambda^{(N)}_{2ij}$
be a Poisson random measure on
${\mathbb R}^+ \times{\mathbb R}^+$ with intensity measure $\frac1 N
l \times l$.
And, for each $i,j\in{\mathbb Z}$, let $\Lambda^{(N)}_{3ij}$ be a
Poisson process with intensity $\frac1 N$.
Then $X_i$ satisfies the following jump equation:
%
%
\begin{eqnarray} \label{eq:construct}
X_i(t) &=& X_i(0) + \int_0^t \Lambda^{(N)}_{1i,b}(du)
- \int_0^t \Lambda^{(N)}_{1i,d}(du) \nonumber\\[-1pt]
&&{} + \sum_j \int_{[0,t]\times[0,\infty)}
\bigl(X_j(u-)-X_i(u-)\bigr)\nonumber\\[-8pt]\\[-8pt]
&&\hspace*{76pt}{}\times
\mathbf{1}_{\{ \xi\leq s(X_j(u-)-X_i(u-)) \}}
\Lambda^{(N)}_{2ij}(du,d\xi) \nonumber\\[-1pt]
&&{} + \sum_j \int_0^t \bigl(X_j(u-)-X_i(u-)\bigr)
\Lambda^{(N)}_{3ij}(du).\nonumber
\end{eqnarray}
%

In the above,\vspace*{-2pt} jumps of $\Lambda^{(N)}_{3ij}$ give possible times
when the type of individual $i$ is replaced by that of individual $j$
due to the resampling mechanism;
jumps of $\Lambda^{(N)}_{2ij}$ give possible times when the type of individual
$i$ is replaced that of individual $j$ due to the selection mechanism;
and jumps of $\Lambda^{(N)}_{1i,b}$ and $\Lambda^{(N)}_{1i,d}$
give possible times when the type
of individual $i$ increases and decreases by 1 due to the beneficial
and deleterious mutation mechanisms, respectively.

In terms of $P_k$, we have
%
%
\begin{eqnarray} \label{eq:Pk}
P_k(t) &=& P_k(0) + \mu\int_0^t q P_{k-1}(u) - P_k(u)\nonumber\\
&&{} + (1-q) P_{k+1}(u)
\, du \nonumber\\[-8pt]\\[-8pt]
&&{} + s\int_0^t \sum_{l\in{\mathbb Z}} (k-l) P_k(u) P_l(u) \, du
\nonumber\\
&&{} + M^{P,1}_k(t) + M^{P,2}_k(t),\nonumber
\end{eqnarray}
where $M^{P,1}_k$ and $M^{P,2}_k$ are orthogonal martingales,
the first arising from the (compensated) mutation mechanism and the
second from
the resampling and selection mechanisms.

We define the conditional quadratic variation
of an $L^2$-martingale $(M_t)_{t\geq0}$ to be the unique previsible
process $\langle M\rangle(t)$ that makes
$M(t)^2-M(0)^2-\langle M\rangle(t)$ a martingale.
See, for example, Chapters II.6 and III.5 of Protter (\citeyear{Protter03}).
With this notation, following the method of, for example,
Ikeda and Watanabe (\citeyear{IkedaWatanabe81}), Section II.3.9, we obtain
%
%
\begin{eqnarray} \label{eq:PkM}
\langle M^{P,1}_k \rangle(t) &=& \frac\mu N \int_0^t
q P_{k-1}(u) + P_k(u) + (1-q) P_{k+1}(u) \, du, \nonumber\\
\langle M^{P,1}_k, M^{P,1}_{k-1} \rangle(t) &=& -\frac\mu N \int_0^t
q P_{k-1}(u) + (1-q) P_k(u) \, du, \nonumber\\
\langle M^{P,1}_k, M^{P,1}_l \rangle(t) &=& 0\qquad
\mbox{if $|k-l|\ge2$}, \\
\langle M^{P,2}_k \rangle(t) &=& \frac1 N \int_0^t \sum_{l\in
{\mathbb Z}}
(2 + s|k-l|) P_k(u) P_l(u) \, du, \nonumber\\
\langle M^{P,2}_k, M^{P,2}_l \rangle(t) &=& - \frac1 N \int_0^t
(2+s|k-l|) P_k(u) P_l(u) \, du \qquad\mbox{if $k\ne l$}.\nonumber
\end{eqnarray}

With the expressions in (\ref{eq:Pk}) and (\ref{eq:PkM}), we can write
the martingale decomposition of the mean
$m(P(t))=\sum_k k P_k(t)$ in the notation of (\ref{def:cn}) as follows
\begin{eqnarray*}
m(P(t)) &=& m(P(0)) + \mu\int_0^t \sum_k k
[q P_{k-1}(u) - P_k(u)\\
&&\hspace*{99.4pt}{} + (1-q) P_{k+1}(u)] \, du \\
&&{} + s\int_0^t \sum_{k,l\in{\mathbb Z}} k(k-l) P_k(u) P_l(u) \, du
+ M^{P,m}(t) \\
&=& m(P(0)) + \mu(2q-1) t\\
&&{} + s\int_0^t c_2(P(u)) \, du + M^{P,m}(t),
\end{eqnarray*}
where $M^{P,m}$ is a martingale, or in
differential notation,
%
%
\begin{equation}\label{eq:dm}
dm(P) = \bigl(\mu(2q-1) + s c_2(P)\bigr) \, dt + d M^{P,m}.
\end{equation}

\section{A nonrigorous argument}
\label{sec3} In this section, we give a nonrigorous argument that leads
to an asymptotic adaptation rate of roughly $\mathcal{O}(\log N)$, as
long as $q$ is strictly positive and regardless of the selection and
mutation parameters. A~rigorous argument in Section \ref{sec4} will
establish a lower bound of $\mathcal{O}(\log^{1-\delta} N)$ on the
adaptation rate.

Our nonrigorous approach is similar to that of Rouzine, Wakeley and
Coffin (\citeyear{RouzineWC03}).
We assume the ``bulk'' of the wave, that is, at $k$'s not too far away
from the mean fitness, behaves like a deterministic traveling wave
and obtain an approximate expression for the shape of this wave.
More specifically, we obtain a set of equations satisfied by
all central moments of the distribution $P$, which will dictate that
the wave is approximately Gaussian. There is, however,
an infinite family of solutions to these equations, parameterized
by the variance of $P$, which ultimately determines the wave speed.
To determine the correct wave speed for a given parameter set
(i.e., population size, mutation and selection coefficients, and the proportion
of beneficial mutations), we use the essentially stochastic
behavior at the front of the wave to calculate the wave speed.
The answer we obtain from both calculations, that is, using the ``bulk'' and
the front of the wave, must be the same. This constraint will yield
an approximate expression for the adaptation rate.

With all martingale terms in (\ref{eq:Pk}) of order $P/N$,
the effect of noise on $P_k$ can be considered to be quite small
if $P_k$ is much larger than $1/N$. For $k$'s where $P_k$ is in this range,
we have from (\ref{eq:Pk}),
%
%
\begin{eqnarray} \label{approx:Pk}
dP_k &\approx& \biggl[ \mu\bigl(q P_{k-1} - P_k + (1-q) P_{k+1}\bigr)
+ s\sum_{l\in{\mathbb Z}} (k-l) P_k P_l \biggr] \, dt \nonumber\\[-8pt]\\[-8pt]
&=& \bigl[ \mu\bigl(q P_{k-1} - P_k + (1-q) P_{k+1}\bigr)
+ s\bigl(k-m(P)\bigr) P_k \bigr] \, dt.\nonumber
\end{eqnarray}
This is similar to (2) in Rouzine, Wakeley and Coffin
(\citeyear{RouzineWC03}).

If we assume that $\{P_k\}_{k\in{\mathbb Z}}$ evolves according to
this deterministic system, then we can calculate the central moments
via the Laplace transform $\psi(\theta; p) = \sum_k e^{\theta
(k-m(p))} p_k$:
\[
d\psi(\theta) = \sum_k e^{\theta(k-m(P))}\, dP_k
- \sum_k \theta e^{\theta(k-m(P))} P_k \,d m(P).
\]
Furthermore, we can obtain from (\ref{eq:dm})
%
%
\begin{equation}\label{approx:mP}
d m(P) \approx \bigl(\mu(2q-1) + s c_2(P)\bigr) \, dt.
\end{equation}
Therefore,
\begin{eqnarray*}
d \psi(\theta) &\approx&\biggl[ \mu\sum_k e^{\theta(k-m(P))}
\bigl(q P_{k-1} - P_k + (1-q) P_{k+1}\bigr)\\
&&\hspace*{5.62pt}{} + s\sum_k e^{\theta(k-m(P))} \bigl(k-m(P)\bigr) P_k \\
&&\hspace*{14.67pt}{} - \sum_k \theta e^{\theta(k-m(P))} P_k
\bigl(\mu(2q-1) + s c_2(P)\bigr) \biggr] \, dt \\
&=& \bigl[ \psi(\theta) \bigl(\mu\bigl(q e^\theta- 1 + (1-q) e^{-\theta}\bigr)
\\
&&\hspace*{28.6pt}{}
- \theta\bigl(\mu(2q-1) + s c_2(P)\bigr)\bigr)
+ s\psi'(\theta) \bigr] \, dt \\
&=& \bigl[ \psi(\theta) \bigl(\mu
\bigl(q e^\theta- 1 + (1-q) e^{-\theta} - \theta(2q-1)\bigr)
- \theta s c_2(P)\bigr)
+ s\psi'(\theta) \bigr] \, dt.
\end{eqnarray*}
We observe that the term with coefficient $\mu$ is $\mathcal
{O}(\theta^2)$,
thus for small $\theta$, the effect of the mutation mechanism on the
centred wave is relatively small compared to the selection mechanism.
We drop the terms arising from the mutation mechanism to obtain
\[
d \psi(\theta) \approx s [ - \psi(\theta) \theta c_2(P)
+ \psi'(\theta) ] \, dt.
\]
Differentiating this repeatedly and
using the fact that $c_n(P) = \psi^{(n)}(0; P)$ for \mbox{$n\ge2$},
we obtain the following approximate system for the central moments $c_n$:
\[
d c_n(P) \approx s\bigl(c_{n+1}(P) - n c_{n-1} (P) c_2(P)\bigr) \, dt.
\]
If we assume the shape of the wave
to be roughly deterministic and stationary, then setting the
expressions on the right-hand side to zero we see that
the central moments of $P$ satisfies
\[
c_n(P) = \cases{
0, &\quad if $n\ge3$ is odd, \vspace*{2pt}\cr
\dfrac{(2n)!}{2^n n!} c_2(P)^{n/2},
&\quad if $n\ge2$ is even,}
\]
which are the central moments of normal distribution with variance $c_2(P)$.
Hence, $P$ is approximately Gaussian, but the variance is not determined.

We can use this information to guess at the asymptotic variance of the wave,
which will also, through (\ref{approx:mP}) yield an
expression for the asymptotic rate of adaptation.
We follow Section 3 of Yu and Etheridge (\citeyear{YuEtheridge08}) and
assume that
$P$ is approximately Gaussian with mean $m(P)$ and variance $b^2$,
and the ``front'' of the wave is approximately where the level of $P$
falls to $1/N$. If the front of the wave is at $K+m(P)$, then
\[
\frac1 {2 \pi b^2} e^{-K^2/2b^2} = \frac1 N,
\]
hence,
%
%
\begin{equation}\label{approx:Kb}
K\approx b \sqrt{2 \log N}.
\end{equation}
To estimate how long it takes the wave to advance by one, we suppose that
a single individual is born at $K+m(P)$ at time zero and estimate the time
it takes for an individual to be born at $K+m(P)+1$. Let $Z(t)$ be the number
of individuals at site $K+m(P)$ at time $t$. Note that these are the fittest
individuals in the population. According to (\ref{approx:Pk}), until a
beneficial mutation falls on site $K+m(P)$, $Z(t)$ increases exponentially
at rate $s K-\mu$. Ignoring beneficial mutations occurring to type
$K-1+m(P)$, that is,
%
%
\begin{equation} \label{eq:Zt}
Z(t) \approx e^{(s K - \mu) t}.
\end{equation}

As the population at site $K+m(P)$ grows, each individual accumulates
beneficial and deleterious mutations at rates $q\mu$ and $(1-q)\mu$,
respectively. The occurrence of the first beneficial mutation will result
in the advance of the wavefront. Using (\ref{eq:Zt}),
we deduce that the probability that no beneficial mutation
occurs to any individuals with fitness type $K+m(P)$ by time $t$ is
\[
\exp\biggl\{ - q \mu\int_0^t Z(u) \, du \biggr\}
= \exp\biggl\{ - \frac{q \mu}{s K - \mu}
\bigl(e^{(s K - \mu) t} -1\bigr) \biggr\},
\]
which gives a wave speed of $(s K-\mu)/\log(s K-\mu)$.

Now we equate the results of our two calculations for the wave speed.
By~(\ref{approx:mP}), the wave speed is
$\mu(2q-1)+s c_2(P) = \mu(2q-1)+s b^2
\approx\mu(2q-1) + s K^2 / (2\log N)$, using the equality involving
$K$ and $b$ in (\ref{approx:Kb}). This leads to the
following consistency condition:
\[
\frac{s K - \mu}{\log(s K - \mu)}
= \mu(2q-1) + \frac{s K^2}{2 \log N}.
\]
For large $K$, this approximately reduces to
%
%
\begin{equation}\label{eq:nonrig}
K \log(s K) = 2 \log N.
\end{equation}
It is easy to see that $K$ must be smaller than $\log N$
but larger than any fractional power of $\log N$.
In fact, (\ref{eq:nonrig}) is a transcendental equation whose solution
can be
written as $K=\frac1 \sigma W(N^{2\sigma})$, where
$W(z)\dvtx[0,\infty)\to[0,\infty)$ is the inverse
function of $z\mapsto z e^z$. Corless et al. (\citeyear{Corlessetal96})
calls the
function $W$ the Lambert $W$ function, and gives useful asymptotic
expansion results of this function near 0 and $\infty$,
for example, (4.20) of Corless et al. (\citeyear{Corlessetal96}).
In particular, the two leading terms of this expansion are
\[
W(z) = \log z - \log\log z + \cdots,
\]
which shows that $K=2\log N - \log(2\log N) + \cdots$ and the
leading term of the wave speed is $2\sigma\log N / (\log\log N)$.
Our rigorous results in Section \ref{sec4} will show that
the rate of adaptation is asymptotically greater than any fractional
power of
$\log N$ as $N\to\infty$.

There are two critical components in the nonrigorous argument
that we presented in this section:
(i) the Gaussian shape of the wave when $N$ is very large, and
(ii) the relation between the speed of the mean and
that of the front of the wave.
The second component above has a rigorous counterpart in
Proposition \ref{thm:mean_tail}, but we have found it difficult to
give a rigorous statement of the shape of the wave
that we can prove and use, therefore our rigorous arguments in Section
\ref{sec4}
does not rely on the first component of the nonrigorous argument.
What takes its place is a comparison argument between the selected
process and the neutral process with only the mutation and resampling
mechanisms.

\section{Stationary measure of the centred process}
\label{sec4}
If $2q-1>0$, then the distribution $P$ tends to move to the right by mutation
and the selection mechanism also works to increase the mean fitness, therefore
no stationary measure for $P$ can exist. If $2q-1<0$, the mutation mechanism
works to decrease the mean fitness but it is
not at all clear the selection mechanism can keep the effects of deleterious
mutations in check and maintain a ``mutation--selection balance.''
However, the process centred about its mean does have a stationary
measure and our first task in this section is to establish this.
Define
\[
\hat p_k = p_{k+m(p)}
\]
for $p\in S^{(N)}$ and $k\in{\mathbb Z}/N$, so that $m(\hat p)=0$ for all
$p$. Define
\begin{eqnarray*}
\hat S^{(N)} &=& \bigl\{ \hat p\dvtx
\mbox{there is some $p \in S^{(N)}$ such that} \\
&&\hspace*{28.32pt} \mbox{$\hat p_k = p_{k+m(p)}$ for all $k\in
{\mathbb Z}/N$} \bigr\}.
\end{eqnarray*}

We observe that every $\hat p\in\hat S^{(N)}$ has all its mass on
points spaced 1 apart and furthermore, the centred process $\hat P$ is
irreducible, that is, all states in $\hat S^{(N)}$ communicate.
To get from state any $\hat p_1$ to any $\hat p_2$, it suffices
to first get to a state where all individuals have the same fitness type.
For example, the following event ensures that at time $t+h$,
all individuals will have the same number of mutations
as carried by individual 1 at time $t$:
%
%
\begin{eqnarray} \label{event:recur}
\max_{i,j} \Lambda^{(N)}_{1i,b}(t,t+h] + \Lambda^{(N)}_{1i,d}(t,t+h]
+ \Lambda^{(N)}_{2ij}(t,t+h] &=& 0 \nonumber\\
\min_i \Lambda^{(N)}_{3i1}(t,t+h] &\ge& 1,\\
\max_i \sum_{j>1} \Lambda^{(N)}_{3ij}(t,t+h] &=& 0.\nonumber
\end{eqnarray}
Then one can get to any configuration in $\hat S^{(N)}$ by
the mutation mechanism alone.
The fact that the event in (\ref{event:recur}) has positive probability
also ensures that the centred process is positive recurrent.
By standard results, for example, Theorem 3.5.3 of Norris (\citeyear{norris98}),
the centred process $\hat P$ is ergodic.
\begin{proposition}
The centred process $(\hat P,X-m(P))$ is ergodic,
that is, there is a unique stationary measure $\pi$ and
regardless of initial condition, the chain converges to the
stationary measure as $t\to\infty$.
\end{proposition}

From now on, we take
\begin{eqnarray*}
\bar S^{(N)} &=& \bigl\{ \bar p\dvtx \mbox{there is some
$p \in S^{(N)}$ and $l\in{\mathbb Z}/N$ such that} \\
&&\hspace*{102pt} \mbox{$\bar p_k = p_{k+l}$
for all $k\in{\mathbb Z}/N$}
\bigr\}
\end{eqnarray*}
to be as our state space for the process $P$ because we may wish to
start the process with an initial configuration that has all its mass
spaced 1 apart but not necessarily falling onto ${\mathbb Z}$. Let
${\mathbb E}^\pi$ denote the expectation started from the stationary
measure~$\pi$. Let $T(t)$ be the semigroup associated with the process
$(P,X)$, then since
\[
\int{\mathbb E}^p [c_2(P(u))] \, d\pi(p) = \int T(u) c_2(p) \, d\pi(p)
= \int c_2(p) \, d\pi(p),
\]
we have
%
%
\begin{eqnarray}\label{eq:Epim}
{\mathbb E}^\pi[m(P(t))]
&=& \int_0^t \int E^{p} [\mu(2q-1) + s c_2(P(u))]
\, d\pi(p) \, du \nonumber\\
&=& \mu(2q-1) t + s\int_0^t \int E^{p} [c_2(\hat P(u))]
\, d\pi(p) \, du \\
&=& \bigl(\mu(2q-1) + s{\mathbb E}^\pi[c_2] \bigr) t.\nonumber
\end{eqnarray}
Thus, it suffices to estimate ${\mathbb E}^\pi[c_2]$ in order to
get a handle on the asymptotic speed at which $m(P)$ increases.
Such an approach resembles the one taken by Desai and Fisher (\citeyear
{DesaiFisher07}).
However, we have found it difficult to estimate ${\mathbb E}^\pi[c_2]$.
Instead, we use the relation between the speed of the mean and
that of the front of the wave. For that, we define
%
%
\begin{equation} \label{def:kcp}
k_c(p)=\max\bigl\{k\dvtx N p_{[k,\infty)} > \log^2 N \bigr\},
\end{equation}
which we view as the location of the front of the wave.
Since $k_c(p)-m(p)=k_c(\hat p)$, we arrive at the following.
\begin{proposition}\label{thm:mean_tail}
For all $t\geq0$, with the stationary measure $\pi$ of the centred process
as the initial condition for the noncentred process $P$, we have
\[
{\mathbb E}^\pi[k_c(P(t))-k_c(P(0))] = {\mathbb E}^\pi[m(P(t))].
\]
\end{proposition}

Roughly speaking, the above proposition states that
the speed of front of the wave is exactly the same as
that of its mean, which seems to be obvious if the wave is of fixed
shape. In the present setting, however, the shape of the wave is stochastic
and this equality holds under the stationary measure of the centred wave.
The idea of relating behavior of the mean and the front of the wave
has been used in our nonrigorous argument in Section \ref{sec3},
as well as in Rouzine, Wakeley and Coffin (\citeyear{RouzineWC03}).
The idea of our main theorem, Theorem \ref{thm:main} below, is to
start the process $P(t)$ from the stationary measure of the centred
process and obtain an $\mathcal{O}(\log^{1-\delta} N)$ lower bound
for the mean
fitness of the population by time 1, as long as the proportion of
beneficial mutations $q$ is strictly positive. In this case, for large
enough population sizes, the mean fitness of the population will
increase at a rate roughly proportional to $\log N$.
All results in what follows are valid for sufficiently large $N$,
which we may not explicitly state all the time.

We first state three results that are
needed for the proof of Theorem \ref{thm:main} below.
Lemma \ref{lem:kc_tail} gives estimates on how far $k_c(P)$ can retreat
on sets of very small probabilities, while Lemma \ref{lem:case0}
compares the selected process with a neutral process to establish that
if the population starts at time 0 with at least $M$ (whose value range
is specified in Lemma \ref{lem:case0}) individuals with
fitness types $\ge K_0$, then the population is expected to have
at least $\log^2 N$ individuals with fitness types
\[
\ge K_0+1.8\log^{1-\varepsilon} M,
\]
where we observe that $C_\mu(\sqrt{c_2(p) N^3} + N^2) e^{-M
e^{-2(1+\mu
)} /4}$
in the statement of Lemma \ref{lem:case0}
is a very small correction factor.
Hence, $k_c(P(1))$ is
expected to be at least $K_0+\log^{1-\varepsilon} M$.
Finally, Proposition \ref{prp:main} states that for any initial condition
$p\in\bar S^{(N)}$, the front is expected to advance at least
$1.7 \log^{1-5\beta} N$ minus a small correction factor.
\begin{lemma}\label{lem:kc_tail}
Let $A^{(N)} \subset\bar S^{(N)}$.
Let $B^{(N)}\in\mathscr{F}_1$ be an event that satisfies
${\mathbb P}^p(B^{(N)})\le\varepsilon(N)$ for all $p\in A^{(N)}$
where $\varepsilon(N)\to0$ as $N\to\infty$.
If $N$ is sufficiently large, then for any $p\in A^{(N)}$,
\[
{\mathbb E}^p \Bigl[\inf_{t\in[0,1]} \min_{i=1,\ldots,N}
\bigl(X_i(t)-X_i(0)\bigr) \mathbf{1}_{B^{(N)}} \Bigr]
\ge- C_\mu\bigl(\sqrt{c_2(p) N^3} + N^2 \bigr) \varepsilon(N)^{1/2},
\]
where $C_\mu$ is a constant depending only on $\mu$.
In particular, if $N$ is sufficiently large, then for any $p\in A^{(N)}$,
\[
{\mathbb E}^p \Bigl[\inf_{t\in[0,1]}
\bigl(k_c(P(t))-k_c(p)\bigr) \mathbf{1}_{B^{(N)}} \Bigr]
\ge- C_\mu\bigl(\sqrt{c_2(p) N^3} + N^2 \bigr) \varepsilon(N)^{1/2}.
\]
The result still holds if we replace the process $(P,X)$ by
the neutral process $(P^{(Y)},Y)$ defined in (\ref{def:Y}).
\end{lemma}
\begin{lemma}\label{lem:case0}
Let $t_1\in[1/2,1]$, $K_0\in{\mathbb Z}/N$, and $\varepsilon\in(0,1)$
be fixed.
Let $M=M(N)$ be a constant that depends on $N$ such that
\[
\frac{M}{e^{\log^{1-0.9\varepsilon} M} \log^2 N} \to\infty
\]
as $N\to\infty$. If $N$ is sufficiently large, then for any $p\in
\bar S^{(N)}$
with $p_{[K_0,\infty)}\ge M/N$, we have
\[
{\mathbb E}^p \Bigl[\inf_{t\in[t_1,1]} k_c(P(t)) \Bigr]
\ge K_0 - C_\mu\bigl(\sqrt{c_2(p) N^3} + N^2 \bigr)
e^{-M e^{-2(1+\mu)} /4}
+ 1.8\log^{1-\varepsilon} M.
\]
\end{lemma}
\begin{proposition}\label{prp:main}
Let $\mu>0$, $q>0$ and $s>0$ be fixed. If $N$ is sufficiently large,
then for any $\beta>0$ and $p\in\hat S^{(N)}$
\begin{eqnarray*}
&&{\mathbb E}^p [k_c(P(1))-k_c(p)] \\
&&\qquad\ge 1.7 \log^{1-5\beta} N
 - C_\mu\bigl( \sqrt{c_2(p) N^3}
+ {\mathbb E}^p \bigl[ \sqrt{c_2(P(t_0)) N^3} \bigr] + N^2 \bigr)
e^{-1/2 \log^2 N},
\end{eqnarray*}
where $t_0=\frac2 s\log^{-\beta} N$.
\end{proposition}
\begin{theorem}\label{thm:main}
Let $\mu>0$, $q>0$ and $s>0$ be fixed. Then for any $\beta>0$,
\[
{\mathbb E}^\pi[m(P(1))] \ge\log^{1-6\beta} N,
\]
if $N$ is sufficiently large.
\end{theorem}
\begin{pf}
We combine Propositions \ref{thm:mean_tail} and \ref{prp:main} to obtain
\begin{eqnarray*}
{\mathbb E}^\pi[m(P(1))] &=& {\mathbb E}^\pi
[k_c(P(1))-k_c(p)] \\
&\ge& 1.7 \log^{1-5\beta} N\\
&&{} - C_\mu\bigl(
N^{3/2} {\mathbb E}^\pi\bigl[ \sqrt{c_2} + \sqrt{c_2(P(t_0))}
\bigr] + N^2 \bigr) e^{-1/2 \log^2 N} \\
&=& 1.7 \log^{1-5\beta} N - C_\mu\bigl(2N^{3/2} {\mathbb E}^\pi
\bigl[\sqrt{c_2}\bigr]
+ N^2 \bigr) e^{-1/2 \log^2 N}.
\end{eqnarray*}
But from (\ref{eq:Epim}), we have
\[
{\mathbb E}^\pi[m(P(1))] = \mu(2q-1) + s{\mathbb E}^\pi[c_2].
\]
Hence,
\begin{eqnarray*}
&&(s+2C_\mu N^{3/2} e^{-1/2 \log^2 N}) {\mathbb E}^\pi[c_2]
\\
&&\qquad\ge1.7 \log^{1-5\beta} N - \mu(2q-1)
- C_\mu N^2 e^{-1/2 \log^2 N},
\end{eqnarray*}
which implies that
\[
{\mathbb E}^\pi[c_2] \ge\frac{1.6}{s} \log^{1-5\beta} N
\]
for sufficiently large $N$. The desired result follows.
\end{pf}

The rest of this work is devoted to the proof of Proposition \ref{prp:main},
which makes use of Lemmas \ref{lem:kc_tail} and \ref{lem:case0}. We define
\begin{eqnarray*}
L&=& \log^{1-3\beta} N, \\
k_d(p) &=& \max\bigl\{k\dvtx N p_{[k,\infty)} > e^{\log^{1-\beta} N} \bigr\}.
\end{eqnarray*}
The number of individuals beyond $k_d$, $e^{\log^{1-\beta} N}$, is much
larger than the number beyond $k_c$ (which is $\log^2 N$) but nevertheless
is only a tiny proportion of the entire population.
The basic idea for the proof of Proposition \ref{prp:main} is to use
Lemma \ref{lem:case0}, which states that
if there are $M$ individuals with fitness types larger than $K_0$ at
time 0,
then $k_c(P)$ is expected to be beyond $K_0+1.8\log^{1-\varepsilon} M$ at
time 1,
where the value of $\varepsilon$ does not depend on $M$ as long as $M$ is
sufficiently large. We can then divide into 2 cases:
(i) if $k_d(P) \ge k_c(p) - L$ before some small time $t_0$
(event $B_1\cup B_2$ below), and (ii)
if $k_d(P) < k_c(p)-L$ throughout the time interval $[0,t_0]$
(event $B_3\cup B_4$ below).
Under case (i), a simple application of Lemma \ref{lem:case0} implies that
the $e^{\log^{1-\beta} N}$ individuals with fitness types larger
than $k_c(p)-L$ are expected to push $k_c(P)$ to beyond
$k_c(p)-L+2L$ at time 1, hence advancing $k_c(P)$ by at least
$L$.
Under case (ii),
the $\log^2 N$ individuals with fitness types larger than $k_c(p)$
will pick off individuals with fitness types smaller than $k_c(p)-L$
(of which there are at least $N-e^{\log^{1-\beta} N}$) via the selection
mechanism at a very fast rate so that with very high probability by
time $t_0$,
$P_{[k_c(p),\infty)}(t_0)$ will be at least $e^{\log^{1-4\beta} N}$.
Lemma \ref{lem:case0} implies that
these $e^{\log^{1-4\beta} N}$ individuals will then push $k_c(P)$ forward
by at least $e^{\log^{1-6\beta} N}$ by time 1. In either case,
the front of the wave moves forward at a high speed.
\begin{pf*}{Proof of Proposition \protect\ref{prp:main}}
We take $t_0=\frac2 s\log^{-\beta} N$ and define
\begin{eqnarray*}
T_0 &=& \inf\{t\ge0\dvtx k_d(P(t)) \ge k_c(p) - L\}, \\
B_1 &=& \bigl\{ P_{[k_c(p) - L,\infty)}(t_0)>e^{\log^{1-2\beta} N},
T_0\le t_0 \bigr\}, \\
B_2 &=& \bigl\{ P_{[k_c(p) - L,\infty)}(t_0)\le e^{\log^{1-2\beta} N},
T_0\le t_0 \bigr\}, \\
B_3 &=& \bigl\{P_{[k_c(p),\infty)}(t_0)>e^{\log^{1-4\beta} N},T_0>t_0 \bigr\},
\\
B_4 &=& \bigl\{P_{[k_c(p),\infty)}(t_0)\le e^{\log^{1-4\beta} N},T_0>t_0
\bigr\}.
\end{eqnarray*}
We will estimate ${\mathbb E}^p [(k_c(P(1))-k_c(p)) \mathbf{1}_B]$
for $B=B_1\cup B_3$ and $B=B_2\cup B_4$.
For $p\in\hat S^{(N)}$ with $k_d(p)\ge k_c(p)-L$, $T_0=0$.
But for those $p$ with $k_d(p)<k_c(p)-L$, we need to
establish that the number of individuals lying in $[k_c(p),\infty)$ grows
quickly, that is, $B_4$ has small probability.
For that, we construct a set valued process $I$
to be dominated by the set of individuals lying in $[k_c(p),\infty)$,
that is, such that $I(t) \subset\{i\dvtx X_i(t)\in[k_c(p),\infty) \}$
for all $t\le T_0$.
Without any loss of generality, we assume that at time 0 individuals
$\{1,\ldots,\log^2 N\}$ lie in $[k_c(p),\infty)$ and
define $I(0)=\{1,\ldots,\log^2 N\}$.
The mechanisms that drive the population $P$ have the following effect
on $I$:
\begin{enumerate}
\item \textit{Mutation}: if any individual $i\in I$ is hit by a deleterious
mutation event, we delete $i$ from $I$.
\item \textit{Selection}: at a selection event when individual $i\in I$ replaces
individual $j$ lying in $(-\infty,k_d(P)]$ at time $t$, we add $j$
to $I$; at a selection event when individual $i\in I$ is replaced by
individual $j\notin I$ at time $t$ (in which case $X_j > X_i$),
we replace $i\in I$ with $j$.
\item \textit{Resampling}: at a resampling event when individual $i$ replaces
individual $j$ at time $t$ (which happens at rate $\frac1 N$),
if $i\notin I$ and $j\in I$ then we delete $j$ from $I$;
if $i\in I$ and $j\notin I$ then we add $j$ to $I$.
\end{enumerate}
Then for $i\in I(t)$, we have $X_i(t)\in[k_c(p),\infty)$, and
$|I|$ has the following transitions:
\begin{enumerate}
\item \textit{Mutation}: $|I|$ decreases by 1 at rate $\mu(1-q) |I|$.
\item \textit{Selection}: $|I|$ increases by 1 at rate
$\frac s N \sum_{i\in I, j\dvtx X_j \le k_d(P)} (X_i-X_j)^+$.
\item \textit{Resampling}: $|I|$ increases by 1 and decreases by 1
both at rate $|I|\frac{N-|I|} N$.
\end{enumerate}
Prior to $T_0$, we have
\begin{eqnarray*}
\frac s N \sum_{i\in I, j\dvtx X_j \le k_d(P)} (X_i-X_j)^+
&\ge& \frac s N \sum_{i\in I, j\dvtx X_j \le k_d(P)} L
\ge\frac s N |I| (N-e^{\log^{1-\beta} N}) L\\
&\ge& 0.9s|I| L
\end{eqnarray*}
for sufficiently large $N$.

Let $Z$ be an integer valued jump process with initial condition
$Z(0)=\log^2 N$ and the following transitions:
\begin{enumerate}
\item $Z$ increases by 1 at rate $0.9sLZ$,
\item $Z$ decreases by 1 at rate $(\mu+1) Z$,
\end{enumerate}
then $Z$ is dominated by $|I|$ before $T_0$.
By Lemma \ref{lem:rapid}(b), if we take $t_0=\frac2 s\log^{-\beta} N$,
which is $\ge\frac1 {0.9sL-\mu-1}(\log\frac{0.9sL}{\mu+1}
+ \log^{1-4\beta} N)$ for sufficiently large $N$,
then
\begin{eqnarray*}
{\mathbb P}^p\bigl(Z(t_0) \le e^{\log^{1-4\beta} N}\bigr)
&\le& \frac1 {(1-e^{-\log^{1-4\beta} N})^{e^{\log^{1-4\beta} N}}}
\biggl(\frac{4(\mu+1)}{0.9sL} \biggr)^{\log^2 N} \\
&\le& C e^{(\log^2 N) (\log C - \log L)} \\
&\le& C e^{-\log^2 N}.
\end{eqnarray*}
Since $|I|$ dominates $Z$ [i.e., $|I(t)|\ge Z(t)$] and
$I$ is dominated by the set of individuals lying in $[k_c(p),\infty)$
before $T_0$, we have
%
%
\begin{equation} \label{ineq:PB3}
{\mathbb P}^p(B_4) \le C e^{-\log^2 N}
\end{equation}
for all $p\in\hat S^{(N)}$.

Now we turn to the event $B_2$. Without any loss of generality, we
assume at time $T_0$, individuals in $A_0=\{1,2,\ldots,\lceil
e^{\log^{1-\beta} N} \rceil\}$ have fitness $\ge k_c(p)-L$. During the
time period $[T_0,t_0]$, the number of resampling events where
individual $i\in A_0$ gets replaced by another individual is
$\operatorname{Poisson}(\frac{N-1}{N} (t_0-T_0))$, so $Y_i$ remains
untouched by a resampling event during $[0,1]$ with probability $\ge
e^{-1}$. Furthermore, no deleterious mutation event falls on $Y_i$
during $[T_0,t_0]$ with probability $e^{-(1-q)\mu(t_0-T_0)} \ge
e^{-\mu}$. Let
\begin{eqnarray*}
A_1 &=& \{i\in A_0\dvtx
X_i \mbox{ remains untouched by a resampling event} \\
&&\hspace*{32.3pt} \mbox{ or a deleterious mutation event during } [T_0,t_0]\},
\end{eqnarray*}
then $|A_1|$ dominates
$\operatorname{Binomial}(\lceil e^{\log^{1-\beta} N} \rceil,e^{-(1+\mu)})$.
By Lemma \ref{lem:hoeff}(a), if\break $p_{[K_0,\infty)}\ge M/N$, then
%
%
\begin{equation}\label{ineq:PB2}
{\mathbb P}^p (B_2) \le e^{-e^{\log^{1-\beta} N} e^{-2(1+\mu)} /2}.
\end{equation}
Combining this and (\ref{ineq:PB3}), implies
%
%
\begin{equation} \label{ineq:PB24}
{\mathbb P}^p (B_2 \cup B_4) \le C e^{-\log^2 N}.
\end{equation}
Hence, by Lemma \ref{lem:kc_tail}, we have for any $p\in\hat S^{(N)}$,
\[
{\mathbb E}^p \bigl[\bigl(k_c(P(1))-k_c(p)\bigr) \mathbf{1}_{B_2 \cup B_4}\bigr]
\ge- C_\mu\bigl(\sqrt{c_2(p) N^3}
+ N^2 \bigr) e^{-1/2 \log^2 N}.
\]

Finally, we turn to events $B_1$ and $B_3$. Both these two events,
unlike $B_2$ and~$B_4$,
will turn out to make large and positive contribution to
the rate of adaptation, and even though we have no estimates on their
probabilities, we expect neither to tend to 0 as $N\to\infty$.
On $B_1$, there are more than
$N e^{\log^{1-2\beta} N}$ individuals in $[k_c(p)-L,\infty)$ at
time $t_0$.
And on $B_3$, at time $t_0$, there are more than $N e^{\log^{1-4\beta
} N}$
individuals in $[k_c(p),\infty)$, therefore for any $p\in\hat S^{(N)}$,
\begin{eqnarray*}
&&\hspace*{-4pt}{\mathbb E}^p \bigl[\bigl(k_c(P(1))-k_c(p)\bigr) \mathbf{1}_{B_1\cup B_3}\bigr] \\
&&\hspace*{-4pt}\qquad= {\mathbb E}^p \bigl[ \{ {\mathbb E}^p [k_c(P(1))
| \mathscr{F}_{t_0} ] - k_c(p) \} \mathbf{1}_{B_1 \cup B_3}
\bigr] \\
&&\hspace*{-4pt}\qquad= {\mathbb E}^p \bigl[ \bigl\{ {\mathbb E}^{P(t_0)}
\bigl[k_c\bigl(P(1-t_0)\bigr)\bigr] - k_c(p) \bigr\} \mathbf{1}_{B_1 \cup B_3}
\bigr] \\
&&\hspace*{-4pt}\qquad\ge {\mathbb E}^p \biggl[
\biggl( -L- C_\mu\bigl(\sqrt{c_2(P(t_0)) N^3} + N^2 \bigr)
e^{-N e^{-2(1+\mu)}/8}
+ 1.8 \log^{1-\beta} \frac N 2 \biggr) \mathbf{1}_{B_1}
\biggr]
\\
&&\hspace*{-4pt}\qquad\quad{} + {\mathbb E}^p \biggl[
\biggl( - C_\mu\bigl(\sqrt{c_2(P(t_0)) N^3} + N^2 \bigr)
e^{-e^{\log^{1-4\beta} N} e^{-2(1+\mu)}/4}
\\
&&\hspace*{-4pt}\qquad\quad\hspace*{139.22pt}{}
+ 1.8\log^{1-\beta} (e^{\log^{1-4\beta} N}) \biggr)
\mathbf{1}_{B_3} \biggr] \\
&&\hspace*{-4pt}\qquad\ge 1.8 (\log^{1-5\beta} N) {\mathbb P}^p(B_1\cup B_3)
- C_\mu e^{-\log^2 N} {\mathbb E}^p \bigl[
\sqrt{c_2(P(t_0)) N^3} + N^2 \bigr],
\end{eqnarray*}
where in the first $\ge$, we use Lemma \ref{lem:case0} twice,
with $K_0=k_c(p)-L$, $M=N/2$, and $\varepsilon=\beta$ for event $B_1$,
and with $K_0=k_c(p)$, $M=e^{\log^{1-4\beta} N}$,
and $\varepsilon=\beta$ for event~$B_3$.

We combine the two estimates above to obtain that if $N$ is
sufficiently large,
then for any $p\in\hat S^{(N)}$,
\begin{eqnarray*}
&& {\mathbb E}^p [k_c(P(1))-k_c(p)] \\
&&\qquad= {\mathbb E}^p \bigl[\bigl(k_c(P(1))-k_c(p)\bigr) \mathbf{1}_{B_1\cup B_3}\bigr]
+ {\mathbb E}^p \bigl[\bigl(k_c(P(1))-k_c(p)\bigr) \mathbf{1}_{B_2\cup B_4}\bigr] \\
&&\qquad\ge 1.8 (\log^{1-5\beta} N) {\mathbb P}^p\bigl((B_2\cup B_4)^c\bigr)
- C_\mu e^{-\log^2 N} {\mathbb E}^p \bigl[
\sqrt{c_2(P(t_0)) N^3} + N^2 \bigr] \\
&&\qquad\quad{} - C_\mu\bigl(\sqrt{c_2(p) N^3} + N^2 \bigr)
e^{-1/2 \log^2 N} \\
&&\qquad\ge 1.7 (\log^{1-5\beta} N)\\
&&\qquad\quad{} - C_\mu\bigl( \sqrt{c_2(p) N^3}
+ {\mathbb E}^p \bigl[ \sqrt{c_2(P(t_0)) N^3} \bigr] + N^2 \bigr)
e^{-1/2 \log^2 N},
\end{eqnarray*}
where we use (\ref{ineq:PB24}) in the last inequality. Hence, we have
the desired result.
\end{pf*}

\section{Proof of supporting lemmas}
\label{sec5}
The lemmas in this section are needed for the proof of Proposition \ref
{prp:main}.
Lemma \ref{lem:hoeff} gives large deviation estimates for
the binomial and Poisson random variables.
Lemma \ref{lem:rapid} establishes a few results on a birth--death process,
which we will use to show that fit individuals pick off unfit
individuals very quickly via the selection mechanism.
We then prove Lemmas \ref{lem:kc_tail} and \ref{lem:case0}.
\begin{lemma}\label{lem:hoeff}
\textup{(a)} Suppose $Z \sim \operatorname{Binomial}(n,\gamma)$, then
$P(Z\le n \gamma/2) \le e^{ - n \gamma^2/2}$.

\textup{(b)} Suppose $\lambda>0$ is fixed and $Z \sim \operatorname{Poisson}(\lambda)$, then
$P(Z\ge n) \hspace*{-0.2pt}\ge\hspace*{-0.2pt}\frac{e^{-\lambda-1/({2n})}}{\sqrt{2\pi}}
( \hspace*{-0.2pt}\frac\lambda{n^2} \hspace*{-0.2pt})^n$. In particular,
if $\varepsilon>0$ is fixed and
$P(Z\ge2\log^{1-\varepsilon} M) \ge c_{(1)} \exp(-\log^{1-0.9\varepsilon
} M)$
for some constant $c_{(1)}$ and sufficiently large $M$.

\textup{(c)} Suppose $Z \sim \operatorname{Poisson}(N\mu)$, then $P(Z\ge N^2) \le C e^{-N\log N}$.
\end{lemma}
\begin{pf}
(a) We use Hoeffding's inequality [Hoeffding (\citeyear{Hoeffding63})]
to prove this:
\[
\qquad
\begin{tabular}{p{320pt}}
Let $X_1,\ldots,X_n$ be i.i.d. random variables taking
values in $[a,b]$. Let\\
$U=X_1+\cdots+X_n$ and $t>0$, then
$P(U-E[U]\ge nt) \le e^{ - 2 n t^2/(b-a)^2 } $.
\end{tabular}
\]
We regard the binomial random variable $n-Z$
as a sum of $n$ independent $\operatorname{Bernoulli}(1-\gamma)$ random variables, then
\begin{eqnarray*}
P(Z\le n \gamma/2) &=& P\bigl(n-Z \ge n (1- \gamma/2)\bigr) \\
&=& P\bigl((n-Z)-n(1-\gamma) \ge n (1- \gamma/2)-n(1-\gamma)\bigr) \\
&=& P\bigl((n-Z)-n(1-\gamma) \ge n \gamma/2\bigr) \\
&\le& e^{ - n \gamma^2/2}
\end{eqnarray*}
by Hoeffding's inequality.

(b) By Stirling's formula [see, e.g., page 257 of Abramowitz and Stegun
(\citeyear{AbSt1965})],
$k! < \sqrt{2\pi} k^{k+1/2} e^{-k+1/({2k})}$ for any integer
$k$. Therefore, for $n\ge1$,
\begin{eqnarray*}
P(Z\ge n) &=& \sum_{k=n}^\infty\frac{e^{-\lambda} \lambda^k}{k!}
\ge e^{-\lambda} \sum_{k=n}^\infty
\frac{\lambda^k}{\sqrt{2\pi} k^{k+1/2} e^{-k+1/({2k})}}
\\
&\ge& e^{-\lambda+n-1/({2n})}
\frac{\lambda^n}{\sqrt{2\pi} n^{n+1/2}}
\ge\frac{e^{-\lambda-1/({2n})}}{\sqrt{2\pi}}
\frac1 {\sqrt n} \biggl(\frac\lambda n \biggr)^n\\
&\ge&\frac{e^{-\lambda-1/({2n})}}{\sqrt{2\pi}}
\biggl( \frac\lambda{n^2} \biggr)^n.
\end{eqnarray*}
We take $n=2\log^{1-\varepsilon} M$, then for sufficiently large $M$,
\begin{eqnarray*}
\frac{e^{-\lambda-1/({2n})}}{\sqrt{2\pi}}
\biggl( \frac\lambda{n^2} \biggr)^n
&\ge&\frac{c_{(1)}}{(4\lambda^{-1}\log^{2-2\varepsilon} M)^{2\log
^{1-\varepsilon} M}}
\\
&=& \frac{c_{(1)}}{ \exp\{ (2\log^{1-\varepsilon} M)
\log(4\lambda^{-1} \log^{2-2\varepsilon} M) \} } \\
&=& c_{(1)} \exp\bigl\{-(\log^{1-0.9\varepsilon} M) (2\log^{-0.1\varepsilon}
M)\\
&&\hspace*{44.2pt}{}\times
\bigl(\log(4\lambda^{-1}) + (2-2\varepsilon)\log\log M\bigr) \bigr\} \\
&\ge& c_{(1)} \exp(-\log^{1-0.9\varepsilon} M).
\end{eqnarray*}

(c) We take $n=N^2$, then
\[
P(Z=n) = e^{-N\mu} \frac{(N\mu)^n}{n!}
\le e^{-N\mu} \frac C {\sqrt n} \biggl(\frac{N\mu e}{n} \biggr)^n
\le C \biggl(\frac1 N \biggr)^N
= C e^{-N\log N},
\]
where we apply Stirling's formula
$k! > \sqrt{2\pi} k^{k+1/2} e^{-k} > c \sqrt k (k/e)^k$.
Consequently,
\begin{eqnarray*}
P(Z\ge n) &=& e^{-N\mu} \sum_{k=n}^\infty\frac{(N\mu)^k}{k!}
\le e^{-N\mu} \frac{(N\mu)^n}{n!} \sum_{k=0}^\infty
\biggl(\frac{N\mu}{n} \biggr)^k \\
&=& P(Z=n) \frac{n}{n-N\mu}
\le C e^{-N\log N}
\end{eqnarray*}
as required.
\end{pf}
\begin{lemma}\label{lem:rapid}
Let $Z$ be an integer valued jump process with initial condition
$Z(0)=Z_0 >0$
and the following transitions:
\begin{enumerate}
\item $Z$ increases by 1 at rate $a Z$,
\item $Z$ decreases by 1 at rate $b Z$,
\end{enumerate}
where $a,b\ge0$ and $a\ne b$, then:

\textup{(a)} For $x\in[0,1)$,
\[
G(x,t) = E(x^{Z_t}) = \biggl( \frac{b(x-1)-(ax-b)e^{-(a-b)t}}{
a(x-1)-(ax-b)e^{-(a-b)t}} \biggr)^{Z_0}.
\]

\textup{(b)} If\vspace*{2pt} $a\ge b$, $M\ge1$ and $t \ge(\log2 \vee\log(aM/b))/(a-b)$, then
$P(Z(t) \le k) \le\frac1 {(1-1/M)^k} ( \frac{4b}{a} )^{Z_0}$.
\end{lemma}
\begin{pf}
(a) It can be shown that $G(x,t)$ satisfies
\[
\frac{\partial}{\partial t} G(x,t) = (ax-b)(x-1) \frac{\partial
}{\partial x} G(x,t)
\]
and that the given $G(x,t)$ satisfies
this PDE with initial condition $G(x,0)=x^{Z_0}$; see, for example,
Theorem 6.11.10 in Grimmett and Stirzaker (\citeyear{GrimmettStirzaker92}).

(b) We take $x=1-1/M$ and apply Markov's inequality to obtain
\begin{eqnarray*}
P\bigl(Z(t) \le k\bigr) &=& P\bigl(x^{Z(t)} \ge x^k\bigr)
\le\frac{E((1-1/M)^{Z(t)})}{(1-1/M)^k} \\
&=& \frac1 {(1-1/M)^k} \biggl(
\frac{b+(a(M-1)-bM)e^{-(a-b)t}}{
a+(a(M-1)-bM)e^{-(a-b)t}} \biggr)^{Z_0} \\
&\le& \frac1 {(1-1/M)^k} \biggl(
\frac{b+a M e^{-(a-b)t}}{a-a e^{-(a-b)t}} \biggr)^{Z_0},
\end{eqnarray*}
where in the last inequality, we use the assumptions $M\ge1$ and
$a\geq b$ to deduce that
$(aM-bM) e^{-(a-b)t}\ge0$. Since $t \ge(\log2 \vee\log(aM/b))/(a-b)$,
we have $a M e^{-(a-b)t}\le b$ and $a e^{-(a-b)t} \le a/2$.
Therefore,
\[
P\bigl(Z(t) \le k\bigr) \le\frac1 {(1-1/M)^k}
\biggl( \frac{4b}{a} \biggr)^{Z_0}
\]
as required.
\end{pf}

Before we prove Lemmas \ref{lem:kc_tail} and \ref{lem:case0}, we first
construct a process $Y$ consisting of individuals that undergo the
mutation and resampling mechanisms of Section~\ref{sec2}
but not the selection mechanism.
Let $Y_i(t)\in{\mathbb Z}/N$, $i=1,\ldots,N$, denote the number of mutations
present in the $i$th individual in the population, then
%
%
\begin{eqnarray} \label{def:Y}
Y_i(t) &=& Y_i(0) + \int_0^t \Lambda^{(N)}_{1i,b}(du)
- \int_0^t \Lambda^{(N)}_{1i,d}(du) \nonumber\\[-8pt]\\[-8pt]
&&{} + \sum_j \int_0^t \bigl(Y_j(u-)-Y_i(u-)\bigr)
\Lambda^{(N)}_{3ij}(du).\nonumber
\end{eqnarray}
Let $P^{(Y)}(t)$ be the empirical measure formed by the $N$ individuals
of the process~$Y$. Since we use the same Poisson random measures and
Poisson processes to construct $X$ and $Y$, we have $Y_i(t)\leq X_i(t)$
for all $t\geq0$ and $i=1,\ldots,N$, provided $Y_i(0)\le X_i(0)$ for
all $i$ at time 0.
\begin{pf*}{Proof of Lemma \protect\ref{lem:kc_tail}}
We prove the result for the neutral process $(P^{(Y)},Y)$, then
since $X$ dominates $Y$, we have the desired result for $(P,X)$ as well.
Let
\[
U=\inf_{t\in[0,1]} \min_{i=1,\ldots,N} \bigl(Y_i(t)-Y_i(0)\bigr)
\]
and $p\in A^{(N)}$ be the initial configuration of the population. We
only need a crude estimate on ${\mathbb E}^p [U \mathbf{1}_{B^{(N)}}
]$. Let $V_1$ be the total number of mutation events (both deleterious
and beneficial) and $V_2$ be the total number of resampling events that
fall on all individuals during $[0,1]$, then $V_1\sim
\operatorname{Poisson} (N\mu)$ and $V_2 \sim
\operatorname{Poisson}(2N)$. Let
\[
k_w(p)=\max\{k-l\dvtx p_k\ne0, p_l \ne0\}
\]
be the width of the support of $p$. Since the resampling mechanism
does not increase the width of the support, $k_w(P(t)) \le k_w(p)+V_1$
for all $t\in[0,1]$.
The most any individual's fitness can decrease due to a resampling
event at $t$ is $k_w(P(t))$, hence
\[
-U \le V_2 \bigl(k_w(p)+V_1\bigr) + (V_2+1) V_1,
\]
where the first term on the right accounts for the possible decrease in
fitness due to each of the $V_2$ resampling events and the second term
accounts for the possible decrease due to mutation events between resampling
events. Hence by Holder's inequality, for any $p\in A^{(N)}$,
\begin{eqnarray*}
{\mathbb E}^p [ |U| \mathbf{1}_{B^{(N)}} ]
&\le& k_w(p) ( {\mathbb E}^p [ V_2^2 ] )^{1/2}
\bigl( {\mathbb P}^p \bigl(B^{(N)}\bigr) \bigr)^{1/2} \\
&&{} + \bigl( {\mathbb E}^p [ (2V_2+1)^4 ] \bigr)^{1/4}
( {\mathbb E}^p [ V_1^4 ] )^{1/4}
\bigl( {\mathbb P}^p \bigl(B^{(N)}\bigr) \bigr)^{1/2} \\
&\le& C_\mu\bigl(k_w(p) N + N^2 \bigr)
\varepsilon(N)^{1/2}.
\end{eqnarray*}
Since $c_2(p) \ge\frac1 N (k_w(p)/2)^2$ for any $p\in\bar S^{(N)}$,
we have the desired result.
\end{pf*}
\begin{pf*}{Proof of Lemma \protect\ref{lem:case0}}
First, we observe that the requirement
%
%
\begin{equation}\label{asn:case0a}
\frac{M}{e^{\log^{1-0.9\varepsilon} M} \log^2 N} \to\infty
\qquad\mbox{as } N\to\infty,
\end{equation}
implies
%
%
\begin{equation}\label{asn:case0b}
\frac M 2 e^{-(1+\mu)} \ge\log^2 N
\end{equation}
for sufficiently large $N$.
Let $Y$ be the neutral process defined in (\ref{def:Y}).
If $p\in\bar S^{(N)}$ and $p_{[K_0,\infty)}\ge M/N$, then at least $M$
individuals lie in $[K_0,\infty)$. Without any loss of generality,
we assume individuals $1,\ldots,M$ lie
in $[K_0,\infty)$. We take the initial condition
$Y_i(0)=X_i(0)$ for all $i=1,\ldots,N$,
then $Y_i(t)\leq X_i(t)$ for all $t\geq0$ and $i=1,\ldots,N$. Let
\begin{eqnarray*}
A_2 &=& \bigl\{i\in\{1,\ldots,M\}\dvtx
Y_i \mbox{ remains untouched by a resampling event} \\
&&\hspace*{74pt} \mbox{ or a deleterious mutation event during } [0,1]\bigr\},
\end{eqnarray*}
then $A_2$ is measurable with respect to the filtration generated by
$\Lambda^{(N)}_{1i,d}$ and $\Lambda^{(N)}_{3ij}$ during the time
period $[0,1]$
and independent from the filtration generated by $\Lambda^{(N)}_{1i,b}$.
Furthermore, the same argument used for (\ref{ineq:PB2}) implies that
the distribution of $|A_2|$ dominates
$\operatorname{Binomial}(M,e^{-(1+\mu)})$. For $p\in\bar S^{(N)}$, we write
%
%
\begin{eqnarray}\label{eq:kc:case1}\qquad
&& {\mathbb E}^p \Bigl[\inf_{t\in[t_1,1]}
k_c\bigl(P^{(Y)}(t)\bigr) \Bigr]
\nonumber\\
&&\qquad= {\mathbb E}^p \biggl[ \inf_{t\in[t_1,1]} k_c\bigl(P^{(Y)}(t)\bigr) \Big|
|A_2|\ge\frac M 2 e^{-(1+\mu)} \biggr]
{\mathbb P}^p \biggl( |A_2|\ge\frac M 2 e^{-(1+\mu)} \biggr)
\\
&&\qquad\quad{} + {\mathbb E}^p \Bigl[ \inf_{t\in[t_1,1]} k_c\bigl(P^{(Y)}(t)\bigr)
\mathbf{1}_{\{ |A_2|< M e^{-(1+\mu)}/2 \}} \Bigr].\nonumber
\end{eqnarray}
By Lemma \ref{lem:hoeff}(a), if $p_{[K_0,\infty)}\ge M/N$, then
%
%
\begin{equation}\label{ineq:B1}
{\mathbb P}^p \biggl( |A_2| < \frac M 2 e^{-(1+\mu)} \biggr)
\le e^{-M e^{-2(1+\mu)} /2}.
\end{equation}

We first deal with the conditional expectation involving the event
$\{ |A_2|\ge\frac M 2 e^{-(1+\mu)} \}$ in (\ref{eq:kc:case1}).
We observe that for the process $Y$ and individuals in $A_2$,
any change in their fitness is due to the beneficial mutation mechanism
and therefore can only increase in time during $[0,1]$.
The number of beneficial mutations that fall on any individual
during $[0,t_1)$ is distributed $\operatorname{Poisson}(q\mu t_1)$ and
since $t_1\ge1/2$, it dominates $\operatorname{Poisson}(q\mu/2)$. Furthermore, it
depends only on $\Lambda^{(N)}_{1i,b}$, therefore is independent of
the set valued random variable $A_2$. Let $K_1$ be
the number of individuals in $A_2$ that have their fitness types
increase by at least
$2\log^{1-\varepsilon} M$ during $[0,t_1]$. If $K_1>\log^2 N$, then
$\inf_{t\in[t_1,1]} k_c(P^{(Y)}(t))\ge K_0+2\log^{1-\varepsilon} M$.
Lemma~\ref{lem:hoeff}(b) with $\lambda=q\mu/2>0$ implies the following:
conditioning on $|A_2|$, the distribution of $K_1$ dominates
$\operatorname{Binomial}(|A_2|,c_{(1)} \exp(-\log^{1-0.9\varepsilon} M))$ for some
constant~$c_{(1)}$,
and then Lemma \ref{lem:hoeff}(a) to obtain
%
%
\begin{eqnarray} \label{ineq:B1a}
&& {\mathbb P}^p \biggl( \inf_{t\in[t_1,1]} k_c\bigl(P^{(Y)}(t)\bigr)
\ge K_0+2\log^{1-\varepsilon} M
\Big| |A_2|\ge\frac M 2 e^{-(1+\mu)} \biggr)
\nonumber\\
&&\qquad \ge{\mathbb P}^p \biggl( K_1 > \log^2 N
\Big| |A_2|\ge\frac M 2 e^{-(1+\mu)} \biggr)
\nonumber\\
&&\qquad \ge{\mathbb P}^p \biggl( K_1
> \frac{c_{(1)} M} 4 e^{-(1+\mu)} e^{-\log^{1-0.9\varepsilon} M}
\Big| |A_2|\ge\frac M 2 e^{-(1+\mu)} \biggr)
\nonumber\\
&&\qquad \ge1-\exp\biggl(-\frac{c_{(1)}^2 M} 4 e^{-(1+\mu)}
e^{-2\log^{1-0.9\varepsilon} M} \biggr) \\
&&\qquad \ge1-\exp\bigl(- c_{(2)} e^{\log M-2\log^{1-0.9\varepsilon} M}
\bigr)
\nonumber\\
&&\qquad \ge1-\exp\bigl(- c_{(2)} e^{0.9 \log M} \bigr),\nonumber
\end{eqnarray}
where $c_{(2)}$ is a constant and we use (\ref{asn:case0b})
in the second inequality. By (\ref{asn:case0a}),
$\frac M 2 e^{-(1+\mu)} > \log^2 N$ for sufficiently large $N$,
hence $\inf_{t\in[t_1,1]} k_c(P^{(Y)}(t)) \ge K_0$
on the event $\{ |A_2|\ge\frac M 2 e^{-(1+\mu)} \}$.
Therefore, (\ref{ineq:B1a}) implies
\[
{\mathbb E}^p \biggl[ \inf_{t\in[t_1,1]} k_c\bigl(P^{(Y)}(t)\bigr) \Big|
|A_2|\ge\frac M 2 e^{-(1+\mu)} \biggr]
\ge K_0 + 1.9\log^{1-\varepsilon} M.
\]

Now we deal with the expectation in (\ref{eq:kc:case1}) involving the event
$\{ |A_2| < \frac M 2 e^{-(1+\mu)} \}$, which, by (\ref{ineq:B1}),
has probability $\le e^{-M e^{-2(1+\mu)} /2}$ if
$p\in\bar S^{(N)}$ and $p_{[K_0,\infty)}\ge M/N$.
We observe that for such $p$, there are more than $\log^2 N$ individuals
with fitness types $\ge K_0$ at time 0, therefore $k_c(p)\ge K_0$. Hence,
Lemma~\ref{lem:kc_tail} implies
\begin{eqnarray*}
&& {\mathbb E}^p \Bigl[\inf_{t\in[t_1,1]} \bigl(k_c\bigl(P^{(Y)}(t)\bigr)-K_0\bigr)
\mathbf{1}_{ \{ |A_2| <  M/2 e^{-(1+\mu)} \} } \Bigr] \\
&&\qquad \ge- C_\mu\bigl(\sqrt{c_2(p) N^3} + N^2 \bigr)
e^{-M e^{-2(1+\mu)} /4},
\end{eqnarray*}
if $p_{[K_0,\infty)}\ge M/N$.
Plugging the above two estimates along with (\ref{ineq:B1})
into (\ref{eq:kc:case1}) yields for $p$ with $p_{[K_0,\infty)}\ge M/N$,
\begin{eqnarray*}
&& {\mathbb E}^p \Bigl[\inf_{t\in[t_1,1]}
k_c\bigl(P^{(Y)}(t)\bigr) \Bigr]
\\
&&\qquad\ge(K_0 + 1.9\log^{1-\varepsilon} M)
{\mathbb P}^p \biggl( |A_2|\ge\frac M 2 e^{-(1+\mu)} \biggr) \\
&&\qquad\quad{} + {\mathbb E}^p \bigl[ \bigl(k_c\bigl(P^{(Y)}(t_1)\bigr)-K_0\bigr)
\mathbf{1}_{\{ |A_2|< M e^{-(1+\mu)}/2 \} } \bigr]\\
&&\qquad\quad\hspace*{0pt}{} + K_0 {\mathbb P}^p\biggl(|A_2|< \frac M 2 e^{-(1+\mu)}\biggr) \\
&&\qquad\ge K_0 - C_\mu\bigl(\sqrt{c_2(p) N^3} + N^2 \bigr)
e^{-M e^{-2(1+\mu)} /4} \\
&&\qquad\quad{} + (1.9\log^{1-\varepsilon} M)
\inf_{p\in\bar S^{(N)}\dvtx p_{[K_0,\infty)}\ge M/N}
{\mathbb P}^p \biggl( |A_2|\ge\frac M 2 e^{-(1+\mu)} \biggr) \\
&&\qquad\ge K_0 - C_\mu\bigl(\sqrt{c_2(p) N^3} + N^2 \bigr)
e^{-M e^{-2(1+\mu)} /4}
+ 1.8\log^{1-\varepsilon} M.
\end{eqnarray*}
Since $X$ dominates $Y$, we have the desired result.
\end{pf*}

\section*{Acknowledgments}
The authors are grateful to Nick Barton for posing the initial problem,
and to Nick Barton and Jonathan Coe for valuable discussions throughout
this research. In addition, we thank three anonymous referees for
valuable advice and suggestions, one of whom pointed out the existence
of a stationary measure for the centred process, while another one
pointed us to the reference on Lambert $W$ functions.

%

%
\printaddresses


\begin{thebibliography}{34}

\bibitem[\protect\citeauthoryear{Abramowitz and Stegun}{1965}]{AbSt1965}
%
\begin{bbook}[auto:SpringerTagBib|2009-01-14|16:51:27]
\bauthor{\bsnm{Abramowitz},~\bfnm{M.}\binits{M.}} \AND
\bauthor{\bsnm{Stegun},~\bfnm{I.~A.}\binits{I.~A.}}
(\byear{1965}).
\btitle{Handbook of Mathematical Functions}.
\bpublisher{Dover}, \baddress{New York}.%
\end{bbook}%
%
\endbibitem%

\bibitem[\protect\citeauthoryear{Barton}{1995}]{Barton95}
%
\begin{barticle}[auto:SpringerTagBib|2009-01-14|16:51:27]
\bauthor{\bsnm{Barton},~\bfnm{N.~H.}\binits{N.~H.}}
(\byear{1995}).
\btitle{Linkage and the limits to natural selection}.
\bjournal{Genetics}
\bvolume{140}
\bpages{821--841}.
\end{barticle}
%
\endbibitem

\bibitem[\protect\citeauthoryear{Barton and Coe}{2009}]{BartonCoe07}
%
\begin{bmisc}[auto:SpringerTagBib|2009-01-14|16:51:27]
\bauthor{\bsnm{Barton},~\bfnm{N.~H.}\binits{N.~H.}} \AND
\bauthor{\bsnm{Coe},~\bfnm{J.~B.}\binits{J.~B.}}
(\byear{2009}).
\bhowpublished{An upper limit to the rate of adaptation. Preprint}.
\end{bmisc}
%
\endbibitem

\bibitem[\protect\citeauthoryear{Brunet et al.}{2006}]{Brunetal06}
%
\begin{barticle}[mr]
\bauthor{\bsnm{Brunet},~\bfnm{E.}\binits{E.}},
\bauthor{\bsnm{Derrida},~\bfnm{B.}\binits{B.}},
\bauthor{\bsnm{Mueller},~\bfnm{A.~H.}\binits{A.~H.}} \AND
\bauthor{\bsnm{Munier},~\bfnm{S.}\binits{S.}}
(\byear{2006}).
\btitle{Noisy traveling waves: Effect of selection on genealogies}.
\bjournal{Europhys. Lett.}
\bvolume{76}
\bpages{1--7}.
\bid{doi={10.1209/epl/i2006-10224-4}, mr={2299937}}
\end{barticle}
%
\endbibitem

\bibitem[\protect\citeauthoryear{Corless et al.}{1996}]{Corlessetal96}
%
\begin{barticle}[mr]
\bauthor{\bsnm{Corless},~\bfnm{R.~M.}\binits{R.~M.}},
\bauthor{\bsnm{Gonnet},~\bfnm{G.~H.}\binits{G.~H.}},
\bauthor{\bsnm{Hare},~\bfnm{D.~E.~G.}\binits{D.~E.~G.}},
\bauthor{\bsnm{Jeffrey},~\bfnm{D.~J.}\binits{D.~J.}} \AND
\bauthor{\bsnm{Knuth},~\bfnm{D.~E.}\binits{D.~E.}}
(\byear{1996}).
\btitle{On the {L}ambert {$W$} function}.
\bjournal{Adv. Comput. Math.}
\bvolume{5}
\bpages{329--359}.
\bid{doi={10.1007/BF02124750}, mr={1414285}}
\end{barticle}
%
\endbibitem

\bibitem[\protect\citeauthoryear{Cuthbertson, Etheridge and Yu}{2009}]{YuEtCul08}
%
\begin{bmisc}[auto:SpringerTagBib|2009-01-14|16:51:27]
\bauthor{\bsnm{Cuthbertson},~\bfnm{C.}\binits{C.}},
\bauthor{\bsnm{Etheridge},~\bfnm{A.~M.}\binits{A.~M.}} \AND
\bauthor{\bsnm{Yu},~\bfnm{Feng}\binits{F.}}
(\byear{2009}).
\bhowpublished{Fixation probability for competing selective sweeps.
Preprint. Available at}
\href{http://arxiv.org/abs/arXiv:0812.0104}{arXiv:0812.0104}.
\end{bmisc}
%
\endbibitem

\bibitem[\protect\citeauthoryear{Desai and Fisher}{2007}]{DesaiFisher07}
%
\begin{barticle}[auto:SpringerTagBib|2009-01-14|16:51:27]
\bauthor{\bsnm{Desai},~\bfnm{M.~M.}\binits{M.~M.}} \AND
\bauthor{\bsnm{Fisher},~\bfnm{D.~S.}\binits{D.~S.}}
(\byear{2007}).
\btitle{Beneficial mutation--selection balance and the effect of
linkage on
positive selection}.
\bjournal{Genetics}
\bvolume{176}
\bpages{1759--1798}.
\bid{doi={10.1534/genetics.106.067678}, pmid={5980116}}
\end{barticle}
%
\endbibitem


\bibitem[\protect\citeauthoryear{Etheridge, Pfaffelhuber and
Wakolbinger}{2009}]{EPW07}
%
\begin{bincollection}[vtex]
\bauthor{\bsnm{Etheridge},~\bfnm{A.~M.}\binits{A.~M.}},
\bauthor{\bsnm{Pfaffelhuber},~\bfnm{P.}\binits{P.}} \AND
\bauthor{\bsnm{Wakolbinger},~\bfnm{A.}\binits{A.}}
(\byear{2009}).
\btitle{How often does the ratchet click? Facts, heuristics,
asymptotics}.
In \bbooktitle{Trends in Stochastic Analysis}.
\bseries{London Math. Soc. Lecture Note Ser.}
\bvolume{353}.
\bpublisher{Cambridge Univ. Press}, \baddress{New York}.
\end{bincollection}
%
\endbibitem

\bibitem[\protect\citeauthoryear{Felsenstein}{1974}]{Felsenstein74}
%
\begin{barticle}[pbm]
\bauthor{\bsnm{Felsenstein},~\bfnm{J.}\binits{J.}}
(\byear{1974}).
\btitle{The evolutionary advantage of recombination}.
\bjournal{Genetics}
\bvolume{78}
\bpages{737--756}.
\bid{pmid={4448362}, pmcid={PMC1213231}}
\end{barticle}
%
\endbibitem

\bibitem[\protect\citeauthoryear{Fisher}{1930}]{Fisher1930}
%
\begin{bbook}[mr]
\bauthor{\bsnm{Fisher},~\bfnm{R.~A.}\binits{R.~A.}}
(\byear{1930}).
\btitle{The Genetical Theory of Natural Selection}.
\bpublisher{Clarendon Press}, \baddress{Oxford}.
\end{bbook}
%
\endbibitem

\bibitem[\protect\citeauthoryear{Gerrish and Lenski}{1998}]{GerrishLenski98}
%
\begin{barticle}[auto:SpringerTagBib|2009-01-14|16:51:27]
\bauthor{\bsnm{Gerrish},~\bfnm{P.~J.}\binits{P.~J.}} \AND
\bauthor{\bsnm{Lenski},~\bfnm{R.~E.}\binits{R.~E.}}
(\byear{1998}).
\btitle{The fate of competing beneficial mutations in an asexual population}.
\bjournal{Genetica}
\bvolume{102/103}
\bpages{127--144}.
\bid{doi={10.1023/A:1017067816551}, pmid={5980116}}
\end{barticle}
%
\endbibitem

\bibitem[\protect\citeauthoryear{Gillespie}{1991}]{Gillespie91}
%
\begin{bbook}[auto:SpringerTagBib|2009-01-14|16:51:27]
\bauthor{\bsnm{Gillespie},~\bfnm{J.~H.}\binits{J.~H.}}
(\byear{1991}).
\btitle{The Causes of Molecular Evolution}.
\bpublisher{Oxford Univ. Press}, \baddress{Oxford, UK}.
\end{bbook}
%
\endbibitem

\bibitem[\protect\citeauthoryear{Gordo and
Charlesworth}{2000}]{GordoCharlesworth00}
%
\begin{barticle}[pbm]
\bauthor{\bsnm{Gordo},~\bfnm{I.}\binits{I.}} \AND
\bauthor{\bsnm{Charlesworth},~\bfnm{B.}\binits{B.}}
(\byear{2000}).
\btitle{On the speed of Muller's ratchet}.
\bjournal{Genetics}
\bvolume{156}
\bpages{2137--2140}.
\bid{pmid={11187462}, pmcid={PMC1461355}}
\end{barticle}
%
\endbibitem

\bibitem[\protect\citeauthoryear{Grimmett and
Stirzaker}{1992}]{GrimmettStirzaker92}
%
\begin{bbook}[mr]
\bauthor{\bsnm{Grimmett},~\bfnm{G.~R.}\binits{G.~R.}} \AND
\bauthor{\bsnm{Stirzaker},~\bfnm{D.~R.}\binits{D.~R.}}
(\byear{1992}).
\btitle{Probability and Random Processes}, \bedition{2nd} ed.
\bpublisher{Clarendon Press}, \baddress{Oxford}.
\bid{mr={1199812}}
\end{bbook}
%
\endbibitem

\bibitem[\protect\citeauthoryear{Haigh}{1978}]{Haigh78}
%
\begin{barticle}[mr]
\bauthor{\bsnm{Haigh},~\bfnm{John}\binits{J.}}
(\byear{1978}).
\btitle{The accumulation of deleterious genes in a population---{M}uller's
ratchet}.
\bjournal{Theoret. Population Biol.}
\bvolume{14}
\bpages{251--267}.
\bid{mr={514423}}
\end{barticle}
%
\endbibitem

\bibitem[\protect\citeauthoryear{Haldane}{1927}]{Haldane27}
%
\begin{barticle}[auto:SpringerTagBib|2009-01-14|16:51:27]
\bauthor{\bsnm{Haldane},~\bfnm{J.~B.~S.}\binits{J.~B.~S.}}
(\byear{1927}).
\btitle{A mathematical theory of natural and artificial selection,
part v:
Selection and mutation}.
\bjournal{Proceedings of the Cambridge Philosophical Society}
\bvolume{23}
\bpages{834--844}.
\end{barticle}
%
\endbibitem

\bibitem[\protect\citeauthoryear{Hegreness et al.}{2006}]{Hegrenessetal06}
%
\begin{barticle}[auto:SpringerTagBib|2009-01-14|16:51:27]
\bauthor{\bsnm{Hegreness},~\bfnm{M.}\binits{M.}},
\bauthor{\bsnm{Shoresh},~\bfnm{N.}\binits{N.}},
\bauthor{\bsnm{Hartl},~\bfnm{D.}\binits{D.}} \AND
\bauthor{\bsnm{Kishony},~\bfnm{R.}\binits{R.}}
(\byear{2006}).
\btitle{An equivalence principle for the incorporation of favourable mutations
in asexual populations}.
\bjournal{Science}
\bvolume{311}
\bpages{1615--1617}.
\end{barticle}
%
\endbibitem

\bibitem[\protect\citeauthoryear{}{1995}]{HiggsWoodcock95}
%
\begin{barticle}[auto:SpringerTagBib|2009-01-14|16:51:27]
\bauthor{\bsnm{Higgs},~\bfnm{Paul}\binits{P.}} \AND
\bauthor{\bsnm{Woodcock},~\bfnm{Glenn}\binits{G.}}
(\byear{1995}).
\btitle{The accumulation of mutations
in asexual
populations, and the structure of genealogical trees in the presence of
selection}.
\bjournal{J. Math. Biol.}
\bvolume{33}
\bpages{677--702}.
\bid{doi={10.1007/BF00184644}, pmid={5980116}}
\end{barticle}
%
\endbibitem

\bibitem[\protect\citeauthoryear{Hill and Robertson}{1966}]{HillRobertson66}
%
\begin{barticle}[auto:SpringerTagBib|2009-01-14|16:51:27]
\bauthor{\bsnm{Hill},~\bfnm{W.~G.}\binits{W.~G.}} \AND
\bauthor{\bsnm{Robertson},~\bfnm{A.}\binits{A.}}
(\byear{1966}).
\btitle{The effect of linkage on limits to artificial selection}.
\bjournal{Genetics Research}
\bvolume{8}
\bpages{269--294}.
\bid{doi={10.1017/S0016672300010156}, pmid={5980116}}
\end{barticle}
%
\endbibitem

\bibitem[\protect\citeauthoryear{Hoeffding}{1963}]{Hoeffding63}
%
\begin{barticle}[mr]
\bauthor{\bsnm{Hoeffding},~\bfnm{Wassily}\binits{W.}}
(\byear{1963}).
\btitle{Probability inequalities for sums of bounded random variables}.
\bjournal{J. Amer. Statist. Assoc.}
\bvolume{58}
\bpages{13--30}.
\bid{mr={0144363}}
\end{barticle}
%
\endbibitem

\bibitem[\protect\citeauthoryear{Ikeda and Watanabe}{1981}]{IkedaWatanabe81}
%
\begin{bbook}[mr]
\bauthor{\bsnm{Ikeda},~\bfnm{Nobuyuki}\binits{N.}} \AND
\bauthor{\bsnm{Watanabe},~\bfnm{Shinzo}\binits{S.}}
(\byear{1981}).
\btitle{Stochastic Differential Equations and Diffusion Processes},
\bedition{2nd} ed.
\bseries{North-Holland Mathematical Library}
\bvolume{24}.
\bpublisher{North-Holland}, \baddress{Amsterdam}.
\bid{mr={1011252}}
\end{bbook}
%
\endbibitem

\bibitem[\protect\citeauthoryear{Muller}{1964}]{Muller64}
%
\begin{barticle}[auto:SpringerTagBib|2009-01-14|16:51:27]
\bauthor{\bsnm{Muller},~\bfnm{H.~J.}\binits{H.~J.}}
(\byear{1964}).
\btitle{The relation of recombination and mutational advance}.
\bjournal{Mutat. Res.}
\bvolume{106}
\bpages{2--9}.%
\end{barticle}%
%
\endbibitem%

\bibitem[\protect\citeauthoryear{Norris}{1998}]{norris98}
%
\begin{bbook}[mr]
\bauthor{\bsnm{Norris},~\bfnm{J.~R.}\binits{J.~R.}}
(\byear{1998}).
\btitle{Markov Chains}.
\bseries{Cambridge Series in Statistical and Probabilistic Mathematics}
\bvolume{2}.
\bpublisher{Cambridge Univ. Press}, \baddress{Cambridge}.
\bid{mr={1600720}}
\end{bbook}
%
\endbibitem

\bibitem[\protect\citeauthoryear{Novella et al.}{1995}]{Novellaetal95}
%
\begin{barticle}[auto:SpringerTagBib|2009-01-14|16:51:27]
\bauthor{\bsnm{Novella},~\bfnm{I.~S.}\binits{I.~S.}},
\bauthor{\bsnm{Elena},~\bfnm{S.~F.}\binits{S.~F.}},
\bauthor{\bsnm{Moya},~\bfnm{A.}\binits{A.}},
\bauthor{\bsnm{Domingo},~\bfnm{E.}\binits{E.}} \AND
\bauthor{\bsnm{Holland},~\bfnm{J.~J.}\binits{J.~J.}}
(\byear{1995}).
\btitle{Size of genetic bottlenecks leading to virus fitness loss is determined
by mean initial population fitness}.
\bjournal{J. Virol.}
\bvolume{69}
\bpages{2869--2872}.
\end{barticle}
%
\endbibitem

\bibitem[\protect\citeauthoryear{Novella et al.}{1999}]{Novellaetal99}
%
\begin{barticle}[auto:SpringerTagBib|2009-01-14|16:51:27]
\bauthor{\bsnm{Novella},~\bfnm{I.~S.}\binits{I.~S.}},
\bauthor{\bsnm{Elena},~\bfnm{S.~F.}\binits{S.~F.}},
\bauthor{\bsnm{Moya},~\bfnm{A.}\binits{A.}},
\bauthor{\bsnm{Domingo},~\bfnm{E.}\binits{E.}} \AND
\bauthor{\bsnm{Holland},~\bfnm{J.~J.}\binits{J.~J.}}
(\byear{1999}).
\btitle{Exponential fitness gains of rna virus populations are limited by
bottleneck effects}.
\bjournal{J. Virol.}
\bvolume{73}
\bpages{1668--1671}.
\end{barticle}
%
\endbibitem

\bibitem[\protect\citeauthoryear{}{2000}]{Orr00}
%
\begin{barticle}[auto:SpringerTagBib|2009-01-14|16:51:27]
\bauthor{\bsnm{Orr},~\bfnm{H.~Allen}\binits{H.~A.}}
(\byear{2000}).
\btitle{The rate of adaptation in asexuals}.
\bjournal{Genetics}
\bvolume{155}
\bpages{961--968}.
\end{barticle}
%
\endbibitem

\bibitem[\protect\citeauthoryear{}{2003}]{Protter03}
%
\begin{bmisc}[auto:SpringerTagBib|2009-01-14|16:51:27]
\bauthor{\bsnm{Protter},~\bfnm{Philip}\binits{P.}}
(\byear{2003}).
\bhowpublished{\textit{Stochastic Integration and Differential
Equations}. \textit{Applications of Mathematics (New York)} \textbf{21}.
Springer, Berlin}.
\end{bmisc}
%
\endbibitem

\bibitem[\protect\citeauthoryear{Rouzine, Brunet and Wilke}{2008}]{RouzineBW07}
%
\begin{barticle}[vtex]
\bauthor{\bsnm{Rouzine},~\bfnm{I.}\binits{I.}},
\bauthor{\bsnm{Brunet},~\bfnm{E.}\binits{E.}} \AND
\bauthor{\bsnm{Wilke},~\bfnm{C.~O.}\binits{C.~O.}}
(\byear{2008}).
\btitle{The traveling-wave approach to asexual evolution: {M}uller's
ratchet and speed of adaptation}.
\bjournal{Theorectial Population Biology}
\bvolume{73}
\bpages{24--46}.
\end{barticle}
%
\endbibitem

\bibitem[\protect\citeauthoryear{Rouzine, Wakeley and
Coffin}{2003}]{RouzineWC03}
%
\begin{barticle}[auto:SpringerTagBib|2009-01-14|16:51:27]
\bauthor{\bsnm{Rouzine},~\bfnm{I.}\binits{I.}},
\bauthor{\bsnm{Wakeley},~\bfnm{J.}\binits{J.}} \AND
\bauthor{\bsnm{Coffin},~\bfnm{J.~M.}\binits{J.~M.}}
(\byear{2003}).
\btitle{The solitary wave of asexual evolution}.
\bjournal{Proc. Nat. Acad. Sci.}
\bvolume{100}
\bpages{587--592}.
\end{barticle}
%
\endbibitem

\bibitem[\protect\citeauthoryear{Stephan, Chao and
Smale}{1993}]{StephanChaoSmale93}
%
\begin{barticle}[auto:SpringerTagBib|2009-01-14|16:51:27]
\bauthor{\bsnm{Stephan},~\bfnm{W.}\binits{W.}},
\bauthor{\bsnm{Chao},~\bfnm{L.}\binits{L.}} \AND
\bauthor{\bsnm{Smale},~\bfnm{J.}\binits{J.}}
(\byear{1993}).
\btitle{The advance of {M}uller's ratchet in a haploid asexual population:
Approximate solution based on diffusion theory}.
\bjournal{Genet. Res.}
\bvolume{61}
\bpages{225--232}.
%
%
\end{barticle}
%
\endbibitem

\bibitem[\protect\citeauthoryear{Taylor}{1996}]{Taylor96}
%
\begin{bbook}[mr]
\bauthor{\bsnm{Taylor},~\bfnm{Michael~E.}\binits{M.~E.}}
(\byear{1996}).
\btitle{Partial Differential Equations: Basic Theory}.
\bseries{Texts in Applied Mathematics}
\bvolume{23}.
\bpublisher{Springer}, \baddress{New York}.
\bid{mr={1395147}}
\end{bbook}
%
\endbibitem

\bibitem[\protect\citeauthoryear{}{2004}]{Wilke04}
%
\begin{barticle}[auto:SpringerTagBib|2009-01-14|16:51:27]
\bauthor{\bsnm{Wilke},~\bfnm{Claus~O.}\binits{C.~O.}}
(\byear{2004}).
\btitle{The speed of adaptation in large asexual populations}.
\bjournal{Genetics}
\bvolume{167}
\bpages{2045--2054}.
%
%
\end{barticle}
%
\endbibitem

\bibitem[\protect\citeauthoryear{}{2008}]{YuEtheridge08}
%
\begin{bincollection}[vtex]
\bauthor{\bsnm{Yu},~\bfnm{Feng}\binits{F.}} \AND
\bauthor{\bsnm{Etheridge},~\bfnm{A.~M.}\binits{A.~M.}}
(\byear{2008}).
\btitle{Rate of adaptation of large
populations}.
In \bbooktitle{Evolutionary Biology from Concept to Application}.
\bpublisher{Springer}, \baddress{Berlin}.
\end{bincollection}
%
\endbibitem

\end{thebibliography}
\end{document}